\documentclass[12pt,reqno]{amsart}

\usepackage{amscd,amsmath,latexsym, thmtools, amssymb, mathtools, enumerate, appendix, amsthm}

\usepackage{fancyhdr}

\usepackage[top=2.5cm, left=3cm, right=3cm, bottom=2.5cm]{geometry}

\usepackage[dvipsnames]{xcolor}
\usepackage{graphicx, tikz, wrapfig, float, caption, subcaption}

\usepackage{tikz-cd}
\usepackage{tkz-euclide}
\usepackage{multicol}

\definecolor{lnkcol}{rgb}{0,0.25,0.75}

\usepackage[msc-links, alphabetic]{amsrefs}

\usepackage{hyperref}
\hypersetup{
  colorlinks=true,
  linkcolor=lnkcol,
  urlcolor=magenta,
  citecolor=magenta
}

\newtheorem{introthm}{Theorem}

\newtheorem{introcor}[introthm]{Corollary}

\newtheorem{thm}{Theorem}[section]
\newtheorem{cor}[thm]{Corollary}
\newtheorem{prop}[thm]{Proposition}
\newtheorem{lem}[thm]{Lemma}
\newtheorem{defn}[thm]{Definition}
\newtheorem{ex}[thm]{Example}

\numberwithin{figure}{section}
\numberwithin{equation}{section}

\title{Quasiconformal Maps between Bowditch Boundaries of Relatively Hyperbolic Groups}

\author{Rana Sardar}
\address{Indian Institute of Science Education and Research Mohali, India}
\email{ranasardar@iisermohali.ac.in}



\begin{document}

\begin{abstract}
    Classifying finitely generated groups up to quasi-isometry is a central problem in geometric group theory. In the context of hyperbolic and relatively hyperbolic groups, one of the key invariants in this classification is the boundary at infinity. Fr\'{e}d\'{e}ric Paulin proved that two hyperbolic groups are quasi-isometric if and only if their Gromov boundaries are quasiconformally equivalent. In this article, we extend this correspondence to relatively hyperbolic groups via their Bowditch boundaries. We introduce a notion of quasiconformal maps on Bowditch boundaries that coarsely preserve shadows of horoballs relative to boundary points. We prove that any coarsely cusp-preserving quasi-isometry between relatively hyperbolic groups induces such a quasiconformal boundary map. Conversely, we prove that every quasiconformal homeomorphism of Bowditch boundaries that coarsely preserves shadows of horoballs arises from a coarsely cusp-preserving quasi-isometry between the groups. 
\end{abstract}


\fontfamily{ptm}\selectfont
\maketitle



\section{Introduction}

Hyperbolic metric spaces, in the sense of Gromov~\cite{Gro87}, generalize the classical notion of negative curvature to a coarse geometric setting. A geodesic metric space \(X\) is called \emph{hyperbolic} if its geodesic triangles are uniformly slim; that is there exists a constant \(\delta \ge 0\) such that for every geodesic triangle in \(X\), each side of the triangle lies within the closed \(\delta\)-neighbourhood of the union of the other two sides. 
A finitely generated group \(G\) is called \emph{hyperbolic} if its Cayley graph, equipped with the word metric associated to some (equivalently, any) finite generating set, is a hyperbolic metric space. 
For a proper hyperbolic metric space \(X\), the \emph{Gromov boundary} \(\partial X\) provides a compactification of \(X\) by adjoining equivalence classes of geodesic rays and plays a central role in rigidity and classification problems. 

A fundamental theme in geometric group theory is the relationship between quasi-isometries and boundary structures. Any quasi-isometry between hyperbolic groups induces a homeomorphism of their boundaries. However, a boundary homeomorphism alone does not determine the quasi-isometry type: there exist examples, due to Pansu, Bourdon, and Pajot, that show hyperbolic groups with homeomorphic boundaries need not be quasi-isometric; see \cites{Bou1997,BP2002,Pan1989}. 
A significant refinement was obtained by Paulin~\cite{Pau96}, who proved that two hyperbolic groups are quasi-isometric if and only if their boundaries are quasiconformally (or quasi-M\"{o}bius) equivalent. 
Quasiconformal structures on ideal boundaries were introduced by Pansu~\cite{Pan1986} in the setting of simply connected Riemannian manifolds of non-positive curvature and were subsequently extended to hyperbolic groups by Paulin \cite{Pau96}. In later work, Pansu~\cite{Pan1989} reformulated the notion of quasiconformality in terms of \emph{shadows of balls}; see \cite{Pan1989}*{Sections~1 and~4}. In this article, we adopt Pansu's shadow-based definition from \cite{Pan1989} for quasiconformal maps on the Gromov boundaries of hyperbolic spaces.
Let $X$ be a proper $\delta$-hyperbolic metric space with Gromov boundary $\partial X$. For a subset $Z \subset X$, the \emph{shadow} relative to a boundary point $\xi \in \partial X$ is the set of endpoints of geodesic rays starting at $\xi$ and meeting $Z$, denoted $\mathcal{O}_\xi(A)$; see \autoref{def_shadow} and \autoref{fig_shadow}. 
Using shadows, one defines a boundary analogue of an annulus. Fix \(R>2\delta\) and \(\xi\in\partial X\). A \emph{ring} in \(\partial X\) is a pair of shadows of balls
\[
\mathcal{O}_\xi B(x,R) \subset \mathcal{O}_\xi B(x',R),
\]
where \(x,x'\in X\) lie on a common bi-infinite geodesic asymptotic to \(\xi\), with \(x'\) lying closer to \(\xi\) than \(x\); see \autoref{def_rings} and \autoref{fig_ring}. If
\(
d_X(x, x') = \log(t)
\) 
for some \(t\geq 1\), then the above ring is called a \emph{\(t\)-ring}. 
We write such a ring as \(B\subset tB,\) or equivalently as \(\frac{1}{t}B\subset B,\) where \(B=\mathcal{O}{\xi}B(x,R)\) and \(tB=\mathcal{O}{\xi}B(x',R)\) (equivalently, \(\tfrac{1}{t}B=\mathcal{O}{\xi}B(x,R)\) and \(B=\mathcal{O}{\xi}B(x',R))\).
The collection of all such rings is denoted by \(\mathcal{R}(\xi,R)\) and is said to be the \emph{ring structure} on \(\partial X\) associated with \(\xi\) and \(R\).
Now let $X$ and $Y$ be two proper hyperbolic metric spaces. Fix $R>2\delta$ and points $\xi \in \partial X$, $\xi' \in \partial Y$, and let $\mathcal{R}(\xi,R)$ and $\mathcal{R}(\xi',R)$ denote the corresponding ring structures on \(\partial X\) and \(\partial Y\), respectively. Given a continuous function $\eta \colon [1,\infty) \to [1,\infty)$, a homeomorphism
\[
f \colon (\partial X, \mathcal{R}(\xi,R)) \longrightarrow (\partial Y, \mathcal{R}(\xi',R)),
\qquad f(\xi) = \xi',
\]
is called \emph{\(\eta\)-quasiconformal} if, for every \(t\)-ring \(B\subset tB\) in \(\mathcal{R}(\xi,R)\), there exists an \(\eta(t)\)-ring $B' \subset \eta(t)B'$ in $\mathcal{R}(\xi',R)$ such that
\[
B' \subset f(B) \subset f(tB) \subset \eta(t)B',
\]
and the same property holds for $f^{-1}$ as well; see \autoref{def_qc_maps}.

Liu~\cite{Liu2022} extended Paulin's theorem to the setting of Morse boundaries by proving that, for two proper cocompact metric spaces whose Morse boundaries contain at least three points, a homeomorphism between their Morse boundaries is induced by a quasi-isometry if and only if both the homeomorphism and its inverse are bi-H\"older, quasisymmetric, or strongly quasiconformal. 
In contrast, Paulin's theorem fails in the broader setting of arbitrary hyperbolic metric spaces. 
A simple counterexample is as follows: Consider the subspace \(X\subset \mathbb{H}^2\) obtained by taking the union of two disjoint closed horoballs connected by a geodesic segment. Then \(X\) is a proper hyperbolic space whose Gromov boundary consists of exactly two points. Although \(\partial X\) is homeomorphic (indeed, quasiconformally or quasi-M\"{o}bius equivalent) to \(\partial\mathbb{R}\), the spaces \(X\) and \(\mathbb{R}\) are not quasi-isometric.

Our goal is to extend this framework to the \emph{Bowditch boundaries of relatively hyperbolic groups}. The notion of \emph{relatively hyperbolic groups}, generalizing Gromov hyperbolic groups, originated in his foundational work \cite{Gro87} and was initially developed by Farb and Bowditch \cites{Farb,Bowditch}. For an overview of the various equivalent definitions of relative hyperbolicity and the relationships among them, we refer the reader to Hruska's work \cite{Hruska}. Throughout this article, we adopt the definition given by Groves and Manning \cite{GM2008}, which is as follows:
Let \(G\) be a finitely generated group with finite generating set \(S\), and let \(X = \mathrm{Cay}(G,S)\) denote its Cayley graph. Let \(\{H_1,\ldots,H_n\}\) be a finite collection of infinite, finitely generated subgroups of infinite index in \(G\) so that $H_i\cap S$ generates $H_i$ for each $i=1,\dots,n$, and let \(\mathcal{H}_G\) denote the collection of all left cosets of $H_1, H_2,..., H_n$. 
We say that \(G\) is \(\delta\)-hyperbolic relative to \(\{H_1,H_2,\dots,H_n\}\) or that \(G\) is \emph{hyperbolic relative to} \(\mathcal{H}_G\) or that $(G,\mathcal{H}_G)$ is \emph{$\delta$-relatively hyperbolic} if the cusped space $X^h$ obtained by attaching combinatorial horoballs $H^h$ to each $H\in \mathcal{H}_G$ is a $\delta$-hyperbolic metric space; see \autoref{def_combinatorial_horoball}, \autoref{def_rel_hyp}.
The boundary \(\partial X^h\), studied extensively by Bowditch~\cite{Bowditch}, is called the \emph{Bowditch boundary} of the relatively hyperbolic group \((G,\mathcal{H}_G)\). For notational convenience, we occasionally denote it by \(\partial G^h\).
For each coset $H \in \mathcal{H}_G$, its limit set in \(\partial X^h\) consists of a single point. A point $p \in \partial X^h$ is called a \emph{parabolic endpoint} if it is the unique limit point of some coset $H_p \in \mathcal{H}_G$. In this case, we say $H_p$ is a \emph{horosphere-like subset} centered at $p$, and its associated combinatorial horoball $H_p^h$ is a \emph{horoball} centered at $p$. Similar to the hyperbolic case, a quasi-isometry between relatively hyperbolic groups that \emph{coarsely preserves cusps}, see \autoref{def_cusp_preserving}, induces a homeomorphism between their Bowditch boundaries; see \autoref{prop_cusp_prev_qi_implies_homeo}, and the converse is not true; see \cites{PS2024,Radhika2025}. Even, like-wise the hyperbolic group setting, the existence of a quasiconformal equivalence between boundaries does not guarantee that the groups are quasi-isometric; see \autoref{ex_qc_not_imply_qi}.

In recent years, several authors have explored conditions under which boundary homeomorphisms between Bowditch boundaries of relatively hyperbolic groups can be lifted to quasi-isometries between those groups. 
Mackay and Sisto introduced the notion of \emph{shadow-respecting quasisymmetric homeomorphisms} between Bowditch boundaries \cite{Mackay-Sisto}*{Definition~2.6}, and proved that such homeomorphisms correspond precisely to coarsely cusp-preserving quasi-isometries between the groups; see \cite{Mackay-Sisto}*{Corollary~1.3}.
In a similar theme, Hruska and Healy \cite{Healy-Hruska}
independently showed that a coarsely cusp preserving quasi-isometry between relatively hyperbolic groups induces both a quasi-isometry between their cusped spaces and a quasisymmetric map between their Bowditch boundaries. 
More recently, Pal and Sardar introduced the notions of \emph{relative cross-ratios} and \emph{relative quasi-M\"{o}bius} homeomorphisms \cite{PS2024}*{Definitions~3.16 and~3.17}, and proved that a relative quasi-M\"{o}bius homeomorphism between Bowditch boundaries of relatively hyperbolic groups induces a coarsely cusp-preserving quasi-isometry between the groups, and vice-versa; see \cite{PS2024}*{Theorems~1 and~2}.

In this article, we generalize Pansu and Paulin's \cites{Pan1989,Pau96} results to the setting of relatively hyperbolic groups. In the context of relatively hyperbolic groups, the associated cusped spaces are Gromov hyperbolic, which allows us to define and study the notions of shadows, rings, and quasiconformal maps on their Bowditch boundaries, analogous to the hyperbolic setting. Moreover, the shadows of horoballs, which arise naturally from the horoball structure of the cusped space, play a crucial role in constructing a quasi-isometry between the groups from a given quasiconformal homeomorphism of their Bowditch boundaries. 
Let \((G,\mathcal{H}_G)\) be a relatively hyperbolic group with Cayley graph $X$, and let $\partial G^h$ denote its Bowditch boundary. Let $\mathcal{R}(\xi, R)$ be a ring structure on \(\partial X^h\) associated to \(R>2\delta\) and $\xi \in \partial X^h$. For each \(H \in \mathcal{H}_G\) centered at a parabolic endpoint \(a_H \ne \xi\), there exists a \(T_1\)-ring \(B_H \subset T_1 B_H\) in $\mathcal{R}(\xi, R)$ that coarsely approximates the shadow of the horoball \(H^h\); that is,
\[
B_H \subset \mathcal{O}_\xi(H^h) \subset T_1B_H,
\]
see \autoref{lem_shadow_of_horoballs}.
Now let $(G_1,\mathcal{H}_{G_1})$ and $(G_2,\mathcal{H}_{G_2})$ be two relatively hyperbolic groups. Suppose $L \ge 1$ is a constant, and
\[
f : (\partial G^h_1,\mathcal{R}(\xi,R)) \longrightarrow (\partial G^h_2,\mathcal{R}(\xi',R))
\]
is an $\eta$-quasiconformal homeomorphism. As defined in \autoref{def_coarsely_shadows_of_horoballs_preserve}, the map $f$ \emph{coarsely preserves shadows of horoballs relative to $\xi$} if the following hold:
\begin{enumerate}[$i)$]
    \item both $f$ and $f^{-1}$ send parabolic endpoints to parabolic endpoints,
    \item for every $H \in \mathcal{H}_{G_1}$ centered at $a \in \partial X^h \setminus \{\xi\}$ and every $H' \in \mathcal{H}_{G_2}$ centered at $f(a) \in \partial Y^h \setminus \{\xi'\}$, suppose
    \[
        B_H \subset \mathcal{O}_\xi(H^h) \subset T_1 B_H
        \quad \text{and} \quad
        B_{H'} \subset \mathcal{O}_{\xi'}(H'^h) \subset T_1 B_{H'}
    \]
    for some $T_1$-rings $B_H \subset T_1 B_H$ in $\mathcal{R}(\xi,R)$ and  
    $B_{H'} \subset T_1 B_{H'}$ in $\mathcal{R}(\xi',R)$. Then

    \medskip
    
    \begin{itemize}
        \item $\displaystyle \tfrac{1}{L} B_{H'} \subset f(B_H)
        \subset f(\mathcal{O}_\xi(H^h))
        \subset f(T_1 B_H)
        \subset L B_{H'}$,

        \medskip
        
        \item $\displaystyle \tfrac{1}{L} B_H \subset f^{-1}(B_{H'})
        \subset f^{-1}(\mathcal{O}_{\xi'}(H'^h))
        \subset f^{-1}(T_1 B_{H'})
        \subset L B_H$.
    \end{itemize}
\end{enumerate}

\bigskip

Our first main result establishes that every coarsely cusp-preserving quasi-isometry between two relatively hyperbolic groups induces a quasiconformal homeomorphism on their Bowditch boundaries that coarsely preserves the shadow of horoballs.

\begin{introthm} \label{thm_main_theorem_1}
    Given $\delta, \epsilon, K \ge 0$, $\lambda \ge 1$, and $R > 2\delta$, there exist a continuous distortion function $\eta \colon [1,\infty) \to [1,\infty)$ and a constant $L \geq 1$ depending only on $\delta,\epsilon,\lambda,K$ and $R$ such that the following holds. 
    Let $(G_1,\mathcal{H}_{G_1})$ and $(G_2,\mathcal{H}_{G_2})$ be two $\delta$-relatively hyperbolic groups for some $\delta \ge 0$. 
    For $R > 2\delta$ and $(\xi,\xi') \in \partial G_1^h \times \partial G_2^h$, let \(\mathcal{R}(\xi,R)\) and \(\mathcal{R}(\xi',R)\) denote the ring structures on $\partial G_1^h$ and $\partial G_2^h$, respectively.    
    Suppose $\varphi \colon G_1 \to G_2$ is a $K$-coarsely cusp-preserving $(\lambda,\epsilon)$-quasi-isometry. Then for every pair of boundary points $(\xi,\xi') \in \partial G_1^h \times \partial G_2^h$ with $\xi' = \partial\varphi^h(\xi)$, the boundary extension 
    \[
    \partial\varphi^h \colon (\partial G_1^h, \mathcal{R}(\xi,R)) \longrightarrow (\partial G_2^h, \mathcal{R}(\xi',R))
    \]
    is an $\eta$-quasiconformal homeomorphism that $L$-coarsely preserves the shadows of horoballs relative to $\xi$.
    Moreover, the distortion function $\eta$ is of the form \(\eta(t) = Ct^\mu,\) for some constants $C \geq 1$ and $\mu \geq 1$ depending only on $\delta, \lambda, \epsilon, K$, and $R$.
\end{introthm}

Left multiplication by a group element of a finitely generated group induces a cusp-preserving isometry (in particular, a $(1,0)$-quasi-isometry) on its Cayley graph. Therefore, by \autoref{thm_main_theorem_1}, for every $R>2\delta$, there exists a distortion function $\eta \colon [1,\infty) \to [1,\infty)$ depending only on $\delta$ and $R$ such that for any element $g$ of a relatively hyperbolic group $(G,\mathcal{H}_G)$, the left multiplication map $L_g$ induces an $\eta$-quasiconformal map \(\partial L_g^h \colon \partial G^h \longrightarrow \partial G^h.\) Thus, the action of a relatively hyperbolic group on itself induces a uniform quasiconformal action on its Bowditch boundary.

The following theorem is the converse of \autoref{thm_main_theorem_1}, showing that such boundary quasiconformal maps that coarsely preserve shadows of horoballs arise precisely from a coarsely cusp-preserving quasi-isometries between the groups.

\begin{introthm}\label{thm_main_theorem_2}
    Let \((G_1,\mathcal{H}_{G_1})\) and \((G_2,\mathcal{H}_{G_2})\) be two \(\delta\)-relatively hyperbolic groups for some \(\delta\geq 0\), and let \( f:\partial G_1^h\to \partial G_2^h \) be a homeomorphism. For \(R>2\delta\) and \((\xi,\xi')\in \partial G_1^h\times \partial G_2^h\), let \(\mathcal{R}(\xi,R)\) and \(\mathcal{R}(\xi',R)\) denote the corresponding ring structures on \(\partial G_1^h\) and \(\partial G_2^h\), respectively. Suppose that there exist constants \(R>2\delta\) and \(L\geq 1\), and a continuous function \(\eta:[1,\infty)\to [1,\infty)\), such that, for every \(\xi\in \partial G_1^h\) with \(\xi'=f(\xi)\), the map
    \[ 
    f:(\partial G_1^h,\mathcal{R}(\xi,R)) \longrightarrow (\partial G_2^h,\mathcal{R}(\xi',R)) 
    \] 
    is \(\eta\)-quasiconformal and \(L\)-coarsely preserves shadows of horoballs relative to \(\xi\).

    Then there exists a \(K\)-coarsely cusp-preserving \((\lambda,\epsilon)\)-quasi-isometry 
    \[
    \varphi_f:G_1\longrightarrow G_2,
    \]
    where the constants \(\lambda\geq 1\), \(\epsilon\geq 0\), and \(K\geq 0\) depend only on \(\delta\), \(\eta\), \(R\), and \(L\). Moreover, the induced boundary homeomorphism \( \partial \varphi_f^h:\partial G_1^h\to \partial G_2^h\) induced by $\partial\varphi_f^h$ coincides with the original homeomorphism \(f\).
\end{introthm}

In the last assertion of \autoref{thm_main_theorem_2}, the term `induced' refers to the boundary homeomorphism arising from a coarsely cusp-preserving quasi-isometry, as described in \autoref{prop_cusp_prev_qi_implies_homeo}.
Note that the quasiconformality of a boundary homeomorphism between Bowditch boundaries of relatively hyperbolic groups does not depend on the choice of ring structure  \(\mathcal{R}(\xi, R)\) used in its definition; see \autoref{cor_qc_is_independet_of_ring_structure}. 
The same was shown by Pansu \cite{Pan1989}*{Proposition 1.11 and 1.12} in the setting of simply connected manifolds with sectional curvature $-a^{2} \le K \leq -b^{2} < 0$, or admitting a cocompact action by their isometry group. 

In view of \autoref{thm_main_theorem_1}, \autoref{thm_main_theorem_2}, \cite{PS2024}*{Theorem 1, Theorem 2} and \cite{Mackay-Sisto}*{Corollary 1.3}, we conclude that the following equivalences hold:

\begin{introcor}
    Let \((G_1,\mathcal{H}_{G_1})\) and \((G_2,\mathcal{H}_{G_2})\) be relatively hyperbolic groups, and let
    \(
    f: \partial G_1^h \longrightarrow \partial G_2^h
    \)
    be a homeomorphism such that both \(f\) and \(f^{-1}\) preserve parabolic endpoints. Then the following are equivalent:
\begin{enumerate}[$(i)$]
    \item \(f\) is a quasiconformal map that coarsely preserves shadows of horoballs relative to each boundary point.
    \item \(f\) is a relative quasi-M\"obius map.
    \item \(f\) is a shadow-respecting quasisymmetric map.
\end{enumerate}
\end{introcor}
\medskip
\paragraph{\textbf{Outline of the paper.}}
In \autoref{sec_Preliminaries}, we review the necessary background on hyperbolic and relatively hyperbolic groups.
In \autoref{sec_qc_maps}, we introduce shadows, rings, and ring structures, and develop the notion of quasiconformal maps on Bowditch boundaries together with their fundamental properties.
In \autoref{sec_shadows_of_horoballs}, we define quasiconformal maps that coarsely preserve shadows of horoballs.
We then prove our first main result, \autoref{thm_main_theorem_1}, showing that coarsely cusp-preserving quasi-isometries induce quasiconformal maps on Bowditch boundaries.
In \autoref{sec_qc_implies_qi_for_cusped_spaces}, we show that every quasiconformal homeomorphism between the Gromov boundaries of visual hyperbolic spaces induces a quasi-isometry between the spaces themselves.
Finally, in \autoref{sec_qc_implies_qi}, we prove our second main result, \autoref{thm_main_theorem_2}, showing that quasiconformal homeomorphisms of Bowditch boundaries that coarsely preserve shadows of horoballs are precisely those arising from coarsely cusp-preserving quasi-isometries.

\medskip
\paragraph{\textbf{Acknowledgment.}}
The author sincerely thanks Prof. Abhijit Pal for many helpful discussions and for his valuable comments. The author is also grateful to the anonymous referee(s) for their careful reading of the manuscript and for their valuable suggestions and comments, which significantly improved the exposition and clarity of this article.

\section{Preliminaries } \label{sec_Preliminaries} 

\subsection{Hyperbolic Metric Spaces}

In a metric space $(X, d_X)$, a \emph{geodesic} between two points $x$ and $y$, denoted by $[x, y]$, is a continuous map $\gamma : [0, d_X(x, y)] \to X$ such that $\gamma(0) = x$, $\gamma(d_X(x, y)) = y$, and $d_X(\gamma(s),\gamma(t))=|s-t| \text{ for all } s,t\in[0,d_X(x,y)]$. The space $(X, d_X)$ is called a \emph{geodesic metric space} if, for every pair of points $x, y \in X$, there exists a geodesic $[x, y]$ joining them. A \emph{geodesic triangle} in $X$ consists of three points $a, b, c \in X$ (called \emph{vertices}) and three geodesic segments $[a, b]$, $[b, c]$, and $[c, a]$ (called \emph{sides}) joining them. A geodesic triangle with vertices $a, b, c \in X$ is denoted by $\triangle(a, b, c)$.

A geodesic metric space \(X\) is called \(\delta\)-\emph{hyperbolic} if every geodesic triangle in \(X\) is \(\delta\)-slim, that is there exists a constant \(\delta \ge 0\) such that for every geodesic triangle in \(X\), each side of the triangle lies within the closed \(\delta\)-neighbourhood of the union of the other two sides, and \emph{hyperbolic} if it is \(\delta\)-hyperbolic for some \(\delta\geq 0\). A finitely generated group is \emph{hyperbolic} if its Cayley graph with respect to some (equivalently, any) finite generating set is hyperbolic. 

For \(\lambda\geq 1\) and \(\epsilon\geq 0\), a map \(\varphi:X\to Y\) is a \((\lambda,\epsilon)\)-\emph{quasi-isometric embedding} if
\[
\frac{1}{\lambda}d_X(x,y)-\epsilon
\leq d_Y(\varphi(x),\varphi(y))
\leq \lambda d_X(x,y)+\epsilon
\]
for all \(x,y\in X\). If, moreover, \(\varphi(X)\) is \(\epsilon\)-coarsely surjective in \(Y\), that is \(Y \subset N_{\epsilon}(\varphi(X))\), then \(\varphi\) is called a \emph{\((\lambda,\epsilon)\)-quasi-isometry}. A \((\lambda,\epsilon)\)-\emph{quasigeodesic} is a \((\lambda,\epsilon)\)-quasi-isometric embedding of an interval into a metric space. For a comprehensive treatment, we refer the reader to \cites{Gro87,BH99}.

In a hyperbolic metric space, any geodesic and any quasigeodesic with the same endpoints remain within a uniformly bounded neighbourhood of each other. This property is known as the \emph{Morse Lemma}, or the \emph{stability of quasigeodesics}. As a consequence, hyperbolicity is invariant under quasi-isometries, and in particular the notion of a hyperbolic group is independent of the choice of finite generating set.
We now state a precise version of the stability of quasigeodesics, including the dependence of the constants.

\begin{prop}[Stability of Quasigeodesics, {\cite{BH99}*{Chapter III.H, Theorem 1.7}}]
\label{Prop_Stability_of_quasigeodesics}
    For all \(\delta \ge 0\), \(\lambda \ge 1\), and \(\epsilon \ge 0\), there exists a constant \( K_{\ref{Prop_Stability_of_quasigeodesics}} = K_{\ref{Prop_Stability_of_quasigeodesics}}(\delta,\lambda,\epsilon) \ge 0\) such that the following holds. If \(X\) is a \(\delta\)-hyperbolic metric space and \(\gamma\) is a \((\lambda,\epsilon)\)-quasigeodesic in \(X\) with endpoints \(x,y \in X\), then \(\gamma\) and any geodesic segment joining \(x\) and \(y\) lie within the \(K_{\ref{Prop_Stability_of_quasigeodesics}}\)-neighbourhood of each other.
\end{prop}


\subsection{Gromov Boundary of a Hyperbolic Metric Space}

A metric space \((X,d_X)\) is called \emph{proper} if every closed ball in \(X\) is compact. Throughout this subsection, let \((X,d_X)\) be a proper hyperbolic metric
space.

\begin{defn}[Gromov product]
Let \(a,b,w \in X\). The \emph{Gromov product} of \(a\) and \(b\) with respect to
\(w\) is defined by
\[
(a,b)_w \coloneqq \frac{1}{2}\bigl(d_X(a,w) + d_X(b,w) - d_X(a,b)\bigr).
\]
\end{defn}

Two geodesic rays \(\gamma_1,\gamma_2 \colon [0,\infty) \to X\) are said to be \emph{asymptotic} if
\[
\sup_{t \ge 0} d_X\bigl(\gamma_1(t),\gamma_2(t)\bigr) < \infty .
\]
Asymptoticity defines an equivalence relation on the set of geodesic rays in \(X\). The equivalence class of a geodesic ray \(\gamma\) is denoted by \([\gamma]\) or by \(\gamma(\infty)\). 
By stability of quasigeodesics \autoref{Prop_Stability_of_quasigeodesics}, given a quasigeodesic ray (or line) $\alpha$ in $X$, there exists a geodesic ray (respectively, line) $\gamma_\alpha$ such that they lie within a uniformly bounded neighbourhood of each other. We define the equivalence class of \(\alpha\) to be the equivalence class \([\gamma_\alpha]\). Since any two such geodesic representatives are asymptotic, this definition is well-defined.

\begin{defn}[Gromov boundary]
    The \emph{Gromov boundary} of \(X\) is the set
    \[
    \partial X \coloneqq \bigl\{ [\gamma] \; \big| \; \gamma \colon [0,\infty) \to X \text{ is a geodesic ray}\bigr\},
    \]
    where \([\gamma]\) denotes the asymptotic equivalence class of \(\gamma\).
\end{defn}

We fix a basepoint \(w \in X\). For \(a \in \partial X\) and \(r \ge 0\), define
\[
V(a,r) \coloneqq \Bigl\{ \, b \in \partial X \; \Big| \; \liminf_{t \to \infty} \bigl(\gamma_1(t),\gamma_2(t)\bigr)_w \ge r, \text{ for some geodesic rays } \gamma_1,\gamma_2 \text{ with }
\]
\[
[\gamma_1]=a,\; [\gamma_2]=b,\; \text{and } \gamma_1(0)=\gamma_2(0)=w \Bigr\}.
\]

Let \(\overline{X} \coloneqq X \cup \partial X\), and equip \(\overline{X}\) with the topology defined in \cite{BH99}*{Definition~3.5}, in which the family $\{V(a,r) \mid r \ge 0\}$ forms a basis of neighbourhoods at $a \in \partial X$. This topology is independent of the choice of basepoint \(w\); see \cite{BH99}*{Proposition~3.7}. With respect to this topology, the Gromov boundary \(\partial X\) is compact. 
Moreover, the boundary \(\partial X\) is \emph{visible}, meaning that any pair of distinct points in \(\partial X\) can be joined by a bi-infinite geodesic; see \cite{BH99}*{Chapter~III.H, Lemma~3.2}. Any two geodesics in \(X\) with the same endpoints in \(\overline{X}\) lie within a \(2\delta\)-neighbourhood of each other; see \cite{BH99}*{Chapter~III.H, Lemma~3.3}.
For a hyperbolic group \(G\), the \emph{Gromov boundary} \(\partial G\) is defined to be the Gromov boundary of any Cayley graph of \(G\) with respect to a finite generating set. This definition is well-defined up to homeomorphism, since the Gromov boundary is invariant under quasi-isometries, as stated in the following theorem.

\begin{thm}[Quasi-isometry invariance, {\cite{BH99}*{Chapter~III.H, Theorem~3.9}}] \label{thm_qi_implies_homeo_hyp}
    Let \(X\) and \(Y\) be proper \(\delta\)-hyperbolic metric spaces for some \(\delta \ge 0\). If \(\varphi \colon X \to Y\) is a quasi-isometric embedding, then the induced boundary map
    \[
    \partial \varphi \colon \partial X \longrightarrow \partial Y, \qquad [\gamma] \mapsto [\varphi \circ \gamma],
    \]
    is a topological embedding. Moreover, if \(\varphi\) is a quasi-isometry, then \(\partial \varphi\) is a homeomorphism.
\end{thm}

\subsection{Relatively Hyperbolic Groups and Boundaries}

In this subsection, we recall the definition of relatively hyperbolic groups in the sense of Groves and Manning \cite{GM2008}. We also introduce the notions of coarsely cusp-preserving maps and Bowditch boundaries. In addition, we collect several fundamental results that will be used later, including the visual boundedness of horoballs and the fact that coarsely cusp-preserving quasi-isometries between relatively hyperbolic groups induce homeomorphisms between their Bowditch boundaries.

\begin{defn}[Combinatorial horoballs, {\cite{GM2008}}]\label{def_combinatorial_horoball}
Let $H$ be a locally finite graph with vertex set $V(H)$ and all edges of length one. Let $d_H$ denote the associated path metric on $H$. The combinatorial horoball based on $H$, denoted by $\mathcal{H}(H,d_H)$, is the graph defined as follows:
\begin{enumerate}
    \item \emph{Vertices}.
    The vertex set is \(\mathcal{H}^{(0)} = H^{(0)} \times \mathbb{Z}_{\ge 0}\), where \(H^{(0)}\) denotes the vertex set of \(H\).

    \item \emph{Edges.}
    The edge set $\mathcal{H}^{(1)}$ consists of two types of edges:
    \begin{itemize}
        \item Vertical edges: For each $n \in \mathbb{N}\cup \{0\}$ and $x \in V(H)$, there is an edge between $(x, n)$ and $(x, n+1)$. 

        \item Horizontal edges: For each $x, y \in V(H)$,
        \begin{itemize}
            \item if \(e\) is an edge of $H$ joining \(x\) to \(y\), then there is a corresponding edge \(\overline{e}\) connecting \((x,0)\) to \((y,0)\), 

            \item  for each $n \in \mathbb{N}$, if $0 < d_H(x, y) \leq 2^n$, then there is a single edge between $(x, n)$ and $(y, n)$. 
        \end{itemize}
    \end{itemize} 
\end{enumerate}
\end{defn}

Equip the graph \(\mathcal{H}(H,d_H)\) with the path metric \(d_{\mathcal H}\) in
which every edge has unit length. With this metric, \(\mathcal{H}(H,d_H)\) is a
hyperbolic metric space; see \cite{GM2008}. For convenience, we will henceforth
denote this space by \((H^{h}, d_{H^{h}})\) and refer to it simply as a
\emph{horoball}, rather than a combinatorial horoball.

\begin{defn}[Relatively hyperbolic groups]\label{def_rel_hyp}
    Let \(G\) be a finitely generated group with finite generating set \(S\), and let \(X=\mathrm{Cay}(G,S)\) be its Cayley graph. Let \(\{H_1,\ldots,H_n\}\) be a finite collection of infinite, finitely generated subgroups of infinite index in \(G\), such that \(H_i \cap S\) generates \(H_i\) for each \(i=1,\ldots,n\). We denote by \(\mathcal{H}_G\) the collection of all left cosets of \(H_1,\ldots,H_n\).
    For a left coset \(H=gH_i \in \mathcal{H}_G\) and points \(a,b \in H\), let \(d_H(a,b)\) denote the distance between \(a\) and \(b\) measured in the word metric of \(H\).     
    To each coset \(H \in \mathcal{H}_G\), we attach a combinatorial horoball \(H^h\), as in \autoref{def_combinatorial_horoball}. The space \(X^h\) obtained by attaching all these horoballs to \(X\) is called the \emph{cusped space} of \(G\).
    
    We say that \(G\) is \(\delta\)-hyperbolic relative to \(\{H_1,H_2,\dots,H_n\}\) or that \(G\) is \(\delta\)-hyperbolic relative to \(\mathcal{H}_G\), or that \((G,\mathcal{H}_G)\) is a \(\delta\)-\emph{relatively hyperbolic group}, if the cusped space \(X^h\) is a \(\delta\)-hyperbolic metric space. The pair \((G,\mathcal{H}_G)\) is called \emph{relatively hyperbolic} if it is \(\delta\)-relatively hyperbolic for some \(\delta \ge 0\).
    The elements of \(\mathcal{H}_G\) are referred to as \emph{horosphere-like subsets}, or simply \emph{horospheres}, and the conjugates of the subgroups \(H_1,\ldots,H_n\) are called the \emph{parabolic subgroups} of \(G\).
\end{defn}

For any \(x,y\in H\in\mathcal H_G\) with \(d_H(x,y)=2^n\), a geodesic in \(H^h\) joining \(x\) and \(y\) has length \(2n+1\) up to a uniformly bounded additive error; see \cite{GM2008}*{Lemma~5.6, Corollary~5.14}. In particular, each \(H \in \mathcal{H}_G\) is uniformly properly embedded in the Cayley graph \(X\). As a consequence, we obtain the following lemma.

\begin{lem}[See e.g. {\cite{Pal}*{Lemmas~1.2.12, 1.2.19}}] \label{lem_Uniformly_Properly_embedded_Lemma}
    Let $(G,\mathcal{H}_G)$ be a $\delta$-relatively hyperbolic group with a Cayley graph $X$, and let $X^h$ denote the associated cusped space, where $\delta \geq 0$. Then \(X\) is properly embedded in \(X^{h}\); that is, for every \(r \ge 0\) there exists a constant \(K_{\ref{lem_Uniformly_Properly_embedded_Lemma}} = K_{\ref{lem_Uniformly_Properly_embedded_Lemma}}(r) \ge 0\) such that for all \(x,y \in X\), if \(d_{X^{h}}(x,y) \le r \), then \(d_X(x,y) \le K_{\ref{lem_Uniformly_Properly_embedded_Lemma}}\).
\end{lem} 

Note that all groups appearing in the pair $(G,\mathcal{H}_G)$ are finitely generated. Hence, the associated cusped space is locally finite and proper (see \cite{GM2008}*{Remark~3.14}).

\begin{defn}[Bowditch boundary, \cite{Bowditch}]
    Let \((G,\mathcal{H}_G)\) be a relatively hyperbolic group with Cayley graph \(X\). The Gromov boundary \(\partial X^{h}\) of the cusped space \(X^{h}\) is called the \emph{Bowditch boundary} of \(G\) (or of \(X\)), and is denoted by \(\partial X^{h}\). 
\end{defn} 

\begin{lem}
    Let \((G,\mathcal{H}_G)\) be a $\delta$-relatively hyperbolic group with a Cayley graph $X$, where $\delta \geq 0$. For each coset \(H \in \mathcal{H}_G\), the Gromov boundary \(\partial H^{h}\) embeds into \(\partial X^{h}\) as a single point. 
\end{lem}

A point \(a \in \partial X^{h}\) is called a \emph{parabolic endpoint} if it arises as the unique boundary point of a horoball \(H^{h}_a\) corresponding to some \(H_a \in \mathcal{H}_G\).

Let \((G,\mathcal{H}_G)\) be a group as in \autoref{def_rel_hyp} (not necessarily relatively hyperbolic). Let \(X\) and \(X'\) be Cayley graphs of \(G\) corresponding to two finite generating sets. The identity map on \(G\) induces an identity \(X \to X'\) which is \(\lambda\)-bi-Lipschitz for some \(\lambda \geq 1\). Since this map fixes each element of \(\mathcal{H}_G\) pointwise, it extends naturally to an identity map between the vertex sets of the associated cusped spaces \(X^h\) and \((X')^h\). Moreover, this extension preserves the lengths of vertical segments and is a quasi-isometry. Consequently, the property of relative hyperbolicity and the associated Bowditch boundary are independent of the choice of finite generating set. 
Dru\c{t}u \cite{Drutu}*{Theorem~1.2} proved that relative hyperbolicity is a quasi-isometry invariant in the following sense. Let \(\phi \colon G \to G'\) be a quasi-isometry between finitely generated groups, and suppose that \(G\) is hyperbolic relative to a finite collection of subgroups \(\{H_1,\ldots,H_n\}\). Then there exists a finite collection of subgroups \(\{H'_1,\ldots,H'_m\}\) of \(G'\) such that \(G'\) is hyperbolic relative to \(\{H'_1,\ldots,H'_m\}\). Moreover, for each \(i \in \{1,\ldots,m\}\), the subgroup \(H'_i\) is quasi-isometrically embedded in some \(H_j\), where \(j=j(i) \in \{1,\ldots,n\}\).

Farb \cite{Farb} proved that, in a Hadamard manifold with pinched negative curvature, horospheres are visually bounded. An analogous result holds for horospheres in relatively hyperbolic groups.

\begin{prop}[Visual diameter bounded, \cite{Farb}*{Lemma 4.4}, \cite{Pal}*{Lemma 1.2.39}] \label{prop_visually_bounded}
    Let $(G,\mathcal{H}_G)$ be a $\delta$-relatively hyperbolic group with a Cayley graph $X$, and let $X^h$ denote the associated cusped space, where $\delta \geq 0$. There exists a constant \(K_{\ref{prop_visually_bounded}} = K_{\ref{prop_visually_bounded}}(\delta) \ge 0\) such that the following holds. For any \(H \in \mathcal{H}_G\) and any geodesic \(\gamma\) in \(X^{h}\) that does not intersect \(H\), let \(I_{\gamma}\) denote the set of points \(p \in H\) for which there exists \(t\) such that the geodesic segment \([\gamma(t),p]\) intersects \(H\) only at \(\{p\}\). Let \(\operatorname{diam}(I_{\gamma})\) denote the diameter of \(I_{\gamma}\). Then \(\sup \bigl\{\, \operatorname{diam}(I_{\gamma}) \;\big|\; \gamma \text{ is a geodesic in } X^{h} \text{ with } \gamma \cap H = \varnothing \,\bigr\}\), called the \emph{visual size} of \(H\), is at most \(K_{\ref{prop_visually_bounded}}\).
\end{prop}

\begin{defn}[Coarsely cusp-preserving quasi-isometry] \label{def_cusp_preserving}
    Let $(G_1,\mathcal{H}_{G_1})$ and $(G_2,\mathcal{H}_{G_2})$ be two relatively hyperbolic groups, with Cayley graphs $X$ and $Y$, respectively. Given $K \ge 0$, a quasi-isometry $\varphi : X \to Y$ is said to be $K$-coarsely cusp-preserving if the following conditions hold:
    \begin{enumerate}[$(i)$]
        \item For every $H \in \mathcal{H}_{G_1}$, there exists $H' \in \mathcal{H}_{G_2}$ such that the image $\varphi(H)$ is contained in the $K$-neighbourhood of $H'$ in $Y$. 

        \item For every $H' \in \mathcal{H}_{G_2}$, there exists $H \in \mathcal{H}_{G_1}$ such that $\varphi^{-1}(H')$ is contained in the $K$-neighbourhood of $H$ in $X$, where $\varphi^{-1}$ is a quasi-isometry inverse of the quasi-isometry $\varphi$.
    \end{enumerate}
\end{defn}

\begin{ex}
    \((i)\) Consider the free group $\mathbb{F}(a,b)$ and the isomorphism $T : \mathbb{F}(a,b) \to \mathbb{F}(a,b)$ defined on generators by $T(a)=ab^{-1}$ and $T(b)=b$. The group $\mathbb{F}(a,b)$ is hyperbolic relative to the cyclic subgroup $\langle aba^{-1}b^{-1} \rangle$. The map $T$ is a quasi-isometry that preserves the cusp subgroup $\langle aba^{-1}b^{-1} \rangle$.
    
    \((ii)\) The groups $(\mathbb{F}(a,b), \langle a \rangle)$ and $(\mathbb{F}(a,b), \langle aba^{-1}b^{-1} \rangle)$ are both relatively hyperbolic. Although the identity map on $\mathbb{F}(a,b)$ is an isometry, it does not preserve the cusp subgroups $\langle a \rangle$ and $\langle aba^{-1}b^{-1} \rangle$.
\end{ex}

\begin{prop}[\cite{Mackay-Sisto}*{Proposition 5.5, Corollary 1.2}, \cite{Pal}*{ Lemma 1.2.31}]\label{prop_cusp_prev_qi_implies_homeo}
    Given $\lambda \ge 1$ and $\epsilon, K \ge 0$, there exist constants $\lambda'=\lambda'(\epsilon, \lambda, K) \ge 1$ and $\epsilon' = \epsilon'( \epsilon, \lambda, K) \ge 0$ such that the following holds. Let $(G_1,\mathcal{H}_{G_1})$ and $(G_2,\mathcal{H}_{G_2})$ be two relatively hyperbolic groups with Cayley graphs $X$ and $Y$, respectively, and suppose that $\varphi: X \to Y$ is a $K$-coarsely cusp-preserving $(\lambda, \epsilon)$-quasi-isometry. Then:
    \begin{enumerate}[$(i)$]
        \item The map $\varphi$ extends to a $(\lambda', \epsilon')$-quasi-isometry $\varphi^h: X^h \to Y^h$ with the following property. Let $H_1 \in \mathcal{H}_{G_1}$ and suppose that $\varphi(H_1)$ is contained in the $K$-neighbourhood of some  $H_2\in \mathcal{H}_{G_2}$. Then for every $(w,t) \in H_1^h$, we have $\varphi^h(w,t)=\gamma_{\varphi(w)}(t)$, where $\gamma_{\varphi(w)}:[0,\infty)\to Y^h$ is a geodesic ray in $Y^h$ with $\gamma_{\varphi(w)}(0)=\varphi(w)$, and there exists $t_0\le K$ such that $\gamma_{\varphi(w)}\big|_{[t_0,\infty)}$ is a vertical geodesic ray in $H_2^h$ with $\gamma_{\varphi(w)}(t_0)\in H_2$.

        \item The boundary extension map $\partial \varphi^h: \partial X^h \to \partial Y^h$ is a homeomorphism.
    \end{enumerate}
\end{prop}

\begin{ex}\label{ex_qc_not_imply_qi}
    A rich class of relatively hyperbolic groups is provided by the \emph{Bianchi groups} \(\Gamma_D := \mathrm{PSL}(2,\mathcal{O}_D),\) which is a concrete family of non-uniform arithmetic lattices in \(\operatorname{Isom}(\mathbb{H}^3_{\mathbb{R}})=\mathrm{PSL}(2,\mathbb{C}),\) the isometry group of real hyperbolic \(3\)-space. Here \(D<0\) is a square-free integer, and \(\mathcal{O}_D\) denotes the ring of integers of the imaginary quadratic field \(\mathbb{Q}(\sqrt{D})\). 
    The quotient \(\mathbb{H}^3_{\mathbb{R}}/\Gamma_D\) is a finite-volume, non-compact hyperbolic \(3\)-orbifold, whose cusps correspond to the parabolic fixed points of \(\Gamma_D\). By Selberg’s Lemma \cite{DK2018}*{Theorem~3.51}, each Bianchi group \(\Gamma_D\) admits a torsion-free finite-index subgroup \(\Gamma_D'\). The quotient \(\mathbb{H}^3_{\mathbb{R}}/\Gamma_D'\) is therefore a finite-volume, non-compact hyperbolic \(3\)-manifold. Consequently, \(\Gamma_D'\) is relatively hyperbolic with respect to its maximal parabolic subgroups. Moreover, the associated cusped space of \(\Gamma_D'\) is quasi-isometric to \(\mathbb{H}^3_{\mathbb{R}}\), and hence its Bowditch boundary is homeomorphic (indeed, quasiconformally equivalent by \autoref{thm_qi_implies_qc}) to the \(2\)-sphere \(\mathbb{S}^2\). We refer the reader to \cite{MR2003} for a detailed discussion of these and other facts of Bianchi groups. 
    
    By a theorem of R.~Schwartz \cite{Schwartz}*{Corollary~1.3}, for any rank-one Lie group \(G\neq \mathrm{Isom}(\mathbb{H}^2)\), two non-uniform lattices in \(G\) are quasi-isometric if and only if they are commensurable. In particular, let \(\Gamma_{D_1}=\mathrm{PSL}(2,\mathcal{O}_{D_1})\) and \(\Gamma_{D_2}=\mathrm{PSL}(2,\mathcal{O}_{D_2})\) be two non-commensurable Bianchi groups, and let \(\Gamma_{D_1}'\) and \(\Gamma_{D_2}'\) be torsion-free finite-index subgroups. Although their Bowditch boundaries are homeomorphic (indeed, quasiconformally equivalent) to \(\mathbb{S}^2\), Schwartz’s rigidity theorem implies that the groups \(\Gamma_{D_1}'\) and \(\Gamma_{D_2}'\) are not quasi-isometric, since they are not commensurable.
\end{ex}

\subsection{Visual Spaces}

\begin{defn}[Visual spaces, {\cite{BS2000}}]
    Let \(X\) be a proper hyperbolic metric space and let \(p\in \overline{X}=X\cup\partial X\). We say that \(X\) is \emph{\(k\)-visual with respect to \(p\)} if for every \(x\in X\) there exists a geodesic ray \(\gamma\) such that \(d(x,\gamma)\le k\), where \(\gamma\) starts at \(p\) when \(p\in X\), and \(\gamma\) is asymptotic to \(p\) when \(p\in\partial X\).

    The space \(X\) is called \emph{visual} if it is \(k\)-visual with respect to some \(p\in\overline{X}\) and some \(k\ge 0\).
\end{defn}

For example, the Poincar\'e disc \(\mathbb{D}^2\) is a hyperbolic metric space that is \(0\)-visual with respect to every point \(p \in \overline{\mathbb{D}}^2\). In contrast, consider the metric space obtained by starting with the real line \(\mathbb{R}\) and, for each natural number \(n\), attaching a vertical line segment of length \(n\) at the point \(n \in \mathbb{R}\). This space \(X\) is a proper \(0\)-hyperbolic metric space, but it is not visual with respect to any point \(p \in \overline{X}\).

The following lemmas are consequences of the slimness of geodesic triangles and
the visibility property of the Gromov boundary.

\begin{lem}\label{lem_visual_change_basepoint}
    Let \(X\) be a proper \(\delta\)-hyperbolic space that is \(k\)-visual with respect to some boundary point \(\xi\in\partial X\), where $\delta,k \geq 0$. Then there exists a constant $k' = k'(\delta,k) \ge 0$ such that $X$ is $k'$-visual with respect to every point $p \in \overline{X}$.
\end{lem}

\begin{proof}
    Let \(x\in X\). Since \(X\) is \(k\)-visual with respect to \(\xi\), there exists a geodesic \(\gamma=[\xi,a]\) such that \(d(x,\gamma)\le k\). For any \(p\in\overline X\), if $p\neq a,b$, consider a geodesic triangle with vertices \(p\), \(\xi\), and \(a\). By \(\delta\)-slimness, every point of \(\gamma\) lies within uniformly bounded distance of one of the other two sides, each of which is a geodesic ray based at (or asymptotic to) \(p\). Hence, the result follows. 
\end{proof}

\begin{lem}\label{lem_existence_of_geodesic_in_hyp}
    Let $G$ be a $\delta$-hyperbolic group with Cayley graph $X$, and suppose that the Gromov boundary $\partial X$ contains at least two points. Then, for every $p \in \overline{X}$, the space $X$ is $3\delta$-visual with respect to $p$. 
\end{lem} 

\begin{proof}
    Since $\partial X$ contains at least two points, there exists a bi-infinite geodesic $\gamma=(a,b)$ with endpoints $a,b\in\partial X$. Fix a vertex $z\in\gamma$. Since $G$ acts transitively on the vertices of its Cayley graph, for every $x\in X$ there exists $g\in G$ such that $x=gz$. Consequently, $x$ lies on the bi-infinite geodesic $g\gamma$, whose endpoints are $ga,gb\in\partial X$. 
    For $p\in\overline X$, if $p\neq a,b$, consider the geodesic triangle with vertices $p$, $ga$, and $gb$. By $\delta$-slimness, every point of $g\gamma=[ga,gb]$ lies within distance $3\delta$ of one of the other two sides, which are geodesic rays based at (or asymptotic to) $p$. Hence, every $x\in X$ lies within distance $3\delta$ of such a ray, and therefore $X$ is $3\delta$-visual with respect to $p$.
\end{proof}

\begin{lem}\label{lem_existence_of_geodesic_in_rel_hyp}
    Let \((G,\mathcal{H}_G)\) be a \(\delta\)-relatively hyperbolic group with Cayley graph \(X\), and let \(X^h\) be the associated cusped space. Then, for every \(p \in \overline{X^h}\), the space \(X^h\) is \(3\delta\)-visual with respect to \(p\).
\end{lem}

\begin{proof}
    Fix \(p\in \overline{X^h}\). Every vertex \(x\in X^h\) lies on a bi-infinite geodesic \(\gamma\). If \(x\) lies in a horoball, one endpoint of \(\gamma\) may be chosen to be the corresponding parabolic point. 
    Let \(a,b\in \partial X^h\) be the endpoints of \(\gamma\). If $p\neq a,b$, consider the ideal geodesic triangle with vertices \(a,b,p\). The point \(x\) lies on the side \([a,b]\). By the \(3\delta\)-slimness of geodesic ideal triangles, every point of \([a,b]\) lies within distance at most \(3\delta\) of the union of the other two sides. Hence \(x\) lies within distance at most \(3\delta\) of either a geodesic ray \([p,a]\) or a geodesic ray \([p,b]\). 
    Therefore, for every \(x\in X^h\), there exists a geodesic ray with endpoint \(p\) whose image passes within distance \(3\delta\) of \(x\). It follows that \(X^h\) is \(3\delta\)-visual with respect to \(p\).
\end{proof}


\section{Quasiconformal Maps on Boundaries at Infinity} \label{sec_qc_maps}

Let \((X,d_X)\) be a proper hyperbolic metric space, and let \(\partial X\) denote its Gromov boundary. If \(\partial X\) consists of a single point, then all boundary-based structures are trivial. Hence, throughout, we assume that \(\partial X\) contains at least two points.
Also recall that for any distinct points \(a,b \in \partial X\), there exists a bi-infinite geodesic in \(X\) with endpoints \(a\) and \(b\), which we denote by \((a,b)\); see \cite{BH99}*{Chapter III.H, Lemma 3.2}. 

\begin{defn}[Shadows, \cite{Pan1989}]\label{def_shadow}
    Let \((X,d_X)\) be a proper hyperbolic metric space and fix a point \(\xi \in \partial X\). For any subset \(Z \subset X\), the \emph{shadow of \(Z\) relative to \(\xi\)}, denoted by \(\mathcal{O}_\xi(Z)\), is the subset of \(\partial X\) defined by
    \[
    \mathcal{O}_\xi(Z) = \bigl\{\, a \in \partial X \setminus \{\xi\} \;\colon \; \text{there exists a geodesic \((a,\xi)\) intersecting \(Z\)} \bigr\}.
    \]
\end{defn}

\begin{figure}[H]
    \centering
    \begin{tikzpicture}[scale = .9]                    
        \filldraw[gray!70] ({2.5*cos 330},{2.5*sin 330}) arc (60:81.8:13);
        \filldraw[gray!70] ({2.5*cos 330},{2.5*sin 330}) -- ({2.5*cos 171.5},{2.5*sin 171.5}) -- ({2.5*cos 198.3},{2.5*sin 198.3});
        \filldraw[gray!70] ({2.5*cos 171.5},{2.5*sin 171.5}) arc (171.5:198.3:2.5);

        \filldraw[gray!30] ({2.5*cos 330},{2.5*sin 330}) arc (60:75:13);
        \filldraw[gray!30] ({2.5*cos 330},{2.5*sin 330}) -- (-1,0.07) -- (-1.1,-.53);

        \filldraw[white] ({2.5*cos 330},{2.5*sin 330}) arc (60:110:5.5);
        \draw[white,thick] ({2.5*cos 330},{2.5*sin 330})--({2.5*cos 199},{-.1+2.5*sin 199});
        \draw[line width =.5 mm] ({2.5*cos 171.5},{2.5*sin 171.5}) arc (171.5:199:2.5);

        \filldraw[gray!140] (-1.1,-.23) circle (.3cm);

        \draw (-1.1,-.23) circle (.3cm);            


        \draw[dashed] ({2.5*cos 330},{2.5*sin 330}) arc (60:81.8:13);
        \draw[dashed] ({2.5*cos 330},{2.5*sin 330}) arc (60:88.1:10);
        \draw[dashed] ({2.5*cos 330},{2.5*sin 330}) arc (60:94.8:8);
        \draw[dashed] ({2.5*cos 330},{2.5*sin 330}) arc (60:102:6.5);        
        \draw[dashed] ({2.5*cos 330},{2.5*sin 330}) arc (60:109:5.5);

        \draw (0,0) circle (2.5 cm);

        \node at ({.25+2.5*cos 330},{-.2+2.5*sin 330}) {$\xi$};
        \node[white] at (-1.1,-.21) {$Z$};
        \node at ({-.75+2.5*cos 180},{-.2+2.5*sin 180}) {$\mathcal{O}_{\xi}(Z)$};



        
        \node[white] at (3.8,.8) {.};   
    \end{tikzpicture}
    \caption{}
    \label{fig_shadow}
\end{figure}


Consider the Poincar\'e upper-half space \(\mathbb{H}^n\), and let \(\gamma\) be a bi-infinite geodesic in \(\mathbb{H}^n\) asymptotic to a point \(\xi \in \partial \mathbb{H}^n = \mathbb{S}^{n-1}\). For any pair of points \(x,x' \in \gamma\) and any \(R > 0\), the shadows \(\mathcal{O}_\xi(B(x,R))\) and \(\mathcal{O}_\xi(B(x',R))\) are nonempty. Moreover, if \(x'\) lies closer to \(\xi\) than \(x\) along \(\gamma\), then \(\mathcal{O}_\xi(B(x,R)) \subset \mathcal{O}_\xi(B(x',R))\).
The following lemma shows that an analogous nesting property holds, in a coarse sense, for arbitrary proper hyperbolic metric spaces.

\begin{lem}[Monotonicity of shadows along a geodesic] \label{lem_Existence_of_rings}
    Let \(\delta \geq 0\), \(R>2\delta\), and let \((X,d_X)\) be a proper \(\delta\)-hyperbolic metric space. Suppose that \((a,\xi)\) is a bi-infinite geodesic in \(X\) with endpoints \(a,\xi\in\partial X\). If \(x_1,x_2\in (a,\xi)\) are such that \(x_1\) lies between \(a\) and \(x_2\) and \(d_X(x_1,x_2)\geq R+\delta\), then
    \[
    \mathcal{O}_\xi B(x_1,R) \subset \mathcal{O}_\xi B(x_2,R).
    \]
\end{lem}

\begin{figure}[H]
        \centering
        \begin{tikzpicture}
            \draw (-4,0)--(7,0);
            \filldraw (0,0) circle (1pt);
            \filldraw (0,.715) circle (1pt);
            \filldraw (4,0) circle (1pt);

            \draw (0,0) circle (.8cm);
            \draw (4,0) circle (.8cm);

            \draw[bend right=3] (-3,1.2) edge (7,0);

            \node at (0,-.3) {$x_1$};
            \node at (4,-.3) {$x_2$};
            \node at (7.3,0) {$\xi$};
            \node at (-3.5,1.3) {$b$};
            \node at (0,1) {$z$};
        \end{tikzpicture}
        \caption{}
\end{figure}

\begin{proof}
    Let \(b \in \mathcal{O}_\xi B(x_1,R)\). Then there exists a geodesic \((b,\xi)\) intersecting the ball \(B(x_1,R)\) at some point \(z \in (b,\xi)\), so that \(d_X(z,x_1) \le R\). Consider the geodesic triangle with vertices \(x_1, z,\xi\). Since \(x_2\) lies on the side \([x_1,\xi)\) and satisfies \(d_X(x_1,x) \ge R+\delta\), the point \(x_2\) is at distance at least \(\delta\) from the segment \([x_1,z]\). By \(\delta\)-slimness of geodesic triangles, the point \(x_2\) lies within distance at most \(2 \delta\) from \([z,\xi)\). Therefore, the geodesic \((b,\xi)\) intersects \(B(x_2,2\delta)\) and so \(B(x_2,R)\), which implies \(b \in \mathcal{O}_\xi B(x_2,R).\) This proves that
    \[
    \mathcal{O}_\xi B(x_1,R) \subset \mathcal{O}_\xi B(x_2,R).
    \]
\end{proof}

\begin{defn}[Rings, \cite{Pan1989}] \label{def_rings}
    Let \((X,d_X)\) be a proper \(\delta\)-hyperbolic metric space for some \(\delta \ge 0\). We fix a point \(\xi \in \partial X\) and a constant \(R > 2\delta\).
    A \emph{ring} centered at a point \(a \in \partial X \setminus \{\xi\}\) is a pair of shadows \(\mathcal{O}_\xi B(x,R) \subset \mathcal{O}_\xi B(x',R)\), where \(x,x' \in (a,\xi)\) are points on a geodesic from \(a\) to \(\xi\), with \(x'\) lying closer to \(\xi\) than \(x\). 

    \begin{figure}[H]
    \centering
    \begin{tikzpicture}                       
        \filldraw[gray!40] ({2.5*cos 330},{2.5*sin 330}) arc (60:81.8:13);
        \filldraw[gray!40] ({2.5*cos 330},{2.5*sin 330}) -- ({2.5*cos 171.5},{2.5*sin 171.5}) -- ({2.5*cos 199},{2.5*sin 199});
        \filldraw[gray!40] ({2.5*cos 171.5},{2.5*sin 171.5}) arc (171.5:198.3:2.5);

        \filldraw[white] ({2.5*cos 330},{2.5*sin 330}) arc (60:110:5.5);
        \draw[white,thick] ({2.5*cos 330},{2.5*sin 330})--({2.5*cos 199},{-.1+2.5*sin 199});

        \draw[line width =.5 mm] ({2.5*cos 171.8},{2.5*sin 171.8}) arc (171.8:164.8:2.5);
        \draw[line width =.5 mm] ({2.5*cos 198.8},{2.5*sin 198.8}) arc (198.8:208.3:2.5);

        \filldraw[gray] (-.37,-.353) circle (.3cm);
        \filldraw[gray] (-1.1,-.23) circle (.3cm);

        \draw (-.37,-.353) circle (.3cm);
        \draw (-1.1,-.23) circle (.3cm);

        \filldraw (-.37,-.32) circle (1pt);
        \filldraw (-1.1,-.21) circle (1pt);
            
        \draw ({2.5*cos 330},{2.5*sin 330}) arc (60:75:19);
        \draw[dashed] ({2.5*cos 330},{2.5*sin 330}) arc (60:81.8:13);
        \draw ({2.5*cos 330},{2.5*sin 330}) arc (60:94.8:8);
        \draw[dashed] ({2.5*cos 330},{2.5*sin 330}) arc (60:109:5.5);
        \draw ({2.5*cos 330},{2.5*sin 330}) arc (60:118:4.5);

        \draw (0,0) circle (2.5 cm);

        \node at ({.25+2.5*cos 330},{-.2+2.5*sin 330}) {$\xi$};
        \node at ({-.35+2.5*cos 180},{-.2+2.5*sin 180}) {$a$};

        \draw[-stealth] (-.37,-.32) --(.5,.6);
        \draw[-stealth] (-1.1,-.23) --(-.23,.6);

        \node at (-.03,.8) {$x$};
        \node at (.7,.8) {$x'$};  
        
        \node[white] at (2.9,.8) {.};
    \end{tikzpicture}
    \caption{}
    \label{fig_ring}
    \end{figure}

    If, in addition, the centers satisfy \(d_X(x,x') = \log t\) for some \(t \ge 1\), then the ring \(\mathcal{O}_\xi B(x,R) \subset \mathcal{O}_\xi B(x',R)\) is called a \emph{\(t\)-ring}.
\end{defn}

We will often abbreviate a \(t\)-ring by writing \(B \subset tB \quad \text{(or equivalently } \tfrac{1}{t}B \subset B\text{)}\), where
\[
B = \mathcal{O}_\xi B(x,R)
\quad \text{and} \quad
tB = \mathcal{O}_\xi B(x',R),
\]
for points \(x,x' \in X\) lying on the same bi-infinite geodesic with endpoint \(\xi\), with \(x'\) closer to \(\xi\) than \(x\).
Equivalently, the notation \(\tfrac{1}{t}B \subset B\) refers to the same configuration, with
\[
\tfrac{1}{t}B = \mathcal{O}_\xi B(x,R)
\quad \text{and} \quad
B = \mathcal{O}_\xi B(x',R).
\]

\begin{defn}[Ring Structures] \label{def_ring_structure}
    Let \((X,d_X)\) be a proper \(\delta\)-hyperbolic metric space, fix a point \(\xi \in \partial X\), and let \(R > 2\delta\).
    The collection of all rings in \(\partial X\) arising from shadows of balls of radius \(R\) relative to \(\xi\) is denoted by \(\mathcal{R}(\xi,R)\).
    We refer to \(\mathcal{R}(\xi,R)\) as the \emph{ring structure} on \(\partial X\) associated to the basepoint \(\xi\) and the radius \(R\).
\end{defn}

The following proposition is a key technical tool that will be used repeatedly in later sections. Roughly speaking, it asserts that in a visual hyperbolic space, if a subset of the boundary can be sandwiched between two comparable rings, then the corresponding centers in the space must lie at uniformly bounded distance from each other.

\begin{prop}\label{Prop_key_prop_rings}
    Let \(\delta,k \ge 0\), let \(R > 2\delta\), and let \(t_1,t_2 \ge 1\). There exists a constant \(M_{\ref{Prop_key_prop_rings}} = M_{\ref{Prop_key_prop_rings}}(\delta,k,R,t_1,t_2) \ge 0\) such that the following holds. 
    Let \((X,d_X)\) be a proper \(\delta\)-hyperbolic metric space which is \(k\)-visual with respect to a point \(\xi \in \partial X\), and let \(\mathcal{R}(\xi,R)\) denote the ring structure on \(\partial X\) associated to \(\xi\) and \(R\). Suppose that \(x_1,x_2 \in X\) and that the rings
    \[
    \mathcal{O}_\xi B(x_1,R) \subset t_1\,\mathcal{O}_\xi B(x_1,R) \quad \text{and} \quad \mathcal{O}_\xi B(x_2,R) \subset t_2\,\mathcal{O}_\xi B(x_2,R) 
    \] 
    belong to \(\mathcal{R}(\xi,R)\). Assume moreover that these rings are mutually comparable in the sense that 
    \[ 
    \mathcal{O}_\xi B(x_1,R) \subset t_2\,\mathcal{O}_\xi B(x_2,R) \quad \text{and} \quad \mathcal{O}_\xi B(x_2,R) \subset t_1\,\mathcal{O}_\xi B(x_1,R). 
    \] 
    Then the centers \(x_1\) and \(x_2\) lie at uniformly bounded distance, namely \(d_X(x_1,x_2) \le M_{\ref{Prop_key_prop_rings}}\).
\end{prop}

\begin{proof}
    For \(i=1,2\), set \(B_i = \mathcal{O}_{\xi}B(x_i,R)\) and \(t_i B_i = \mathcal{O}_{\xi}B(x_i',R)\). By the definition of a \(t_i\)-ring, the points \(x_i\) and \(x_i'\) lie on a bi-infinite geodesic $\gamma_i$ in \(X\) asymptotic to \(\xi\), that is, \(\gamma_i(-\infty)=\xi\), with \(x_i'\) closer to \(\xi\) and 
    \[ 
    d_X(x_i,x_i') = \log(t_i).
    \]     
    From the given inclusions of rings, we have 
    \[ 
    \mathcal{O}_{\xi}B(x_1,R) \subset \mathcal{O}_{\xi}B(x_2',R) \quad \text{and} \quad \mathcal{O}_{\xi}B(x_2,R) \subset \mathcal{O}_{\xi}B(x_1',R). 
    \] 
    By \(\delta\)-hyperbolicity, this implies that \(\gamma_1\) intersects the \(2\delta\)-neighborhood of \(B(x_2',R)\), and similarly \(\gamma_2\) intersects the \(2\delta\)-neighborhood of \(B(x_1',R)\). Then there are three possible cases: the ball $B(x_1',R + 2\delta)$ intersects $\gamma_2$ either on the segment $[x_2,x_2']$ or on $[x_2,\gamma_2(\infty))$ or on $[x_2',\xi)$.

    \medskip
    \noindent\textbf{Case 1.} Suppose \(B(x_1',R+2\delta)\) intersects \(\gamma_2\) at a point \(q \in [x_2,x_2']\); see \autoref{fig_E(x,t)_is_bounded_1}.  

    \begin{figure}[H]
        \centering
        \begin{tikzpicture}
            \draw[bend left=10] (0,1) edge (9,0);
            \filldraw (1.5,1.05) circle (1pt);
            \filldraw (4.7,.94) circle (1pt);
            \draw (1.5,1.05) circle (.5 cm);
            \draw (4.7,.94) circle (.6 cm);

            
            \draw[bend left=20] (0,-1) edge (9,0);
            \filldraw (1.5,-.34) circle (1pt);
            \filldraw (6.2,.45) circle (1pt);
            \draw (1.5,-.34) circle (.5 cm);
            \draw (6.2,.45) circle (.6 cm);
            \filldraw (4.7,.44) circle (1pt);

            \node at (1.5,1.25) {$x_1$};
            \node at (4.7,1.2) {$x_1'$};
            \node at (1.5,-.6) {$x_2$};
            \node at (6.2,0.2) {$x_2'$};
            \node at (4.8,0.1) {$q$};

            \node at (9.3,0) {$\xi$};   
            \node at (7.7,1.2) {$B(x_2',R + 2\delta)$};
            \node at (0.2,-.2) {$B(x_2,R)$};

            \node at (-.5,1) {$\gamma_1$};
            \node at (-.5,-.9) {$\gamma_2$};
        \end{tikzpicture} 
        \caption{}
        \label{fig_E(x,t)_is_bounded_1}
    \end{figure}   
    
    Then 
    \[ 
    \begin{aligned} 
        d_X(x_1,x_2) & \le d_X(x_1,x_1') + d_X(x_1',q) + d_X(q,x_2) \\ 
        &\le \log(t_1) + (R+2\delta) + \log(t_2). 
    \end{aligned}
    \]  

    \medskip
    \noindent\textbf{Case 2.} Suppose \(B(x_1',R+2\delta)\) intersects \(\gamma_2\) at a point \(q \in [x_2,\gamma_2(\infty))\); see \autoref{fig_E(x,t)_is_bounded_2}. We claim that 
    \[ 
    d_X(q,x_2) \le k + 2R + 6\delta + 1. 
    \]

    \begin{figure}[H]
        \centering        
        \begin{tikzpicture}
            \draw[bend left=10] (-2,.5) edge (9,0);
            \filldraw (-1.3,.58) circle (1pt);
            \filldraw (1.3,.8) circle (1pt);
            \draw (-1.3,.58) circle (.5 cm);
            \draw (1.3,.8) circle (.6 cm);

            
            \draw[bend left=10] (-2.5,-.3) edge (9,0);
            \filldraw (4,.45) circle (1pt);
            \filldraw (6.5,.32) circle (1pt);
            \draw (4,.45) circle (.5 cm);
            \draw (6.5,.32) circle (.6 cm); 
            \filldraw (2.8,.42) circle (1pt);
            \node at (2.8,0.12) {$z$};

            \node at (-1.3,.8) {$x_1$};
            \node at (1.3,1.15) {$x_1'$};
            \node at (4,.2) {$x_2$};
            \node at (6.5,0) {$x_2'$};

            \node at (9.3,0) {$\xi$};       
            \node at (3,1.5) {$B(x_1',R + 2\delta)$};

            \node at (-2.5,.5) {$\gamma_1$};
            \node at (-3,-.3) {$\gamma_2$};

            \filldraw (1.3,0.32) circle (1pt);
            \node at (1.3,-.1) {$q$};
        \end{tikzpicture} 
        \caption{}
        \label{fig_E(x,t)_is_bounded_2}
    \end{figure} 

    Assume, for contradiction, that 
    \[ 
    d_X(q,x_2) > k + 2R + 6\delta + 1. 
    \] 
    Let \(R_0 = k + 2\delta + 1\), and choose a point \(z \in [q,x_2]\) such that \(d_X(z,x_2) = R_0 + \delta\). By \autoref{lem_Existence_of_rings}, 
    \[ 
    \mathcal{O}_{\xi}B(z,R_0) \subset \mathcal{O}_{\xi}B(x_2,R). 
    \]
    Since \(d_X(x_1',q) \le R + 2\delta\) and \(d_X(q,z) > 2R + 3\delta\), again by \autoref{lem_Existence_of_rings} we have 
    \[ 
    \mathcal{O}_{\xi}B(x_1',R) \subset \mathcal{O}_{\xi}B(q,2R+2\delta) \subset \mathcal{O}_{\xi}B(z,2\delta) \subset \mathcal{O}_{\xi}B(z,R_0). 
    \]

    On the other hand, since \(X\) is \(k\)-visual with respect to \(\xi\), there exists a bi-infinite geodesic \((a,\xi)\), with \(a \in \partial X\), which intersects \(B(z,R_0)\) but avoids \(B(z,2\delta)\). Consequently, 
    \[ 
    a \in \mathcal{O}_{\xi}B(z,R_0) \subset \mathcal{O}_{\xi}B(x_2,R), \quad \text{but} \quad a \notin \mathcal{O}_{\xi}B(z,2\delta) \supseteq \mathcal{O}_{\xi}B(x_1',R), 
    \] 
    contradicting the inclusion \(\mathcal{O}_{\xi}B(x_2,R) \subset \mathcal{O}_{\xi}B(x_1',R)\). This proves the claim.

    Therefore, 
    \[ 
    \begin{aligned} 
        d_X(x_1,x_2) & \le d_X(x_1,x_1') + d_X(x_1',q) + d_X(q,x_2) \\ 
        & \le \log(t_1) + (R+2\delta) + (k + 2R + 6\delta + 1). 
    \end{aligned} 
    \] 

    \medskip 
    \noindent\textbf{Case 3.} Suppose \(B(x_1',R+2\delta)\) intersects \(\gamma_2\) on the segment \([x_2',\xi)\). Then \(B(x_2',R+2\delta)\) intersects \(\gamma_1\) either on \([x_1,x_1']\) or on \([x_1,\gamma_1(\infty))\). Interchanging the roles of \((x_1,x_1')\) and \((x_2,x_2')\) in Cases~1 and~2, we obtain 
    \[ 
    d_X(x_1,x_2) \le \log(t_1) + \log(t_2) + (R+2\delta), 
    \] 
    or 
    \[ 
    d_X(x_1,x_2) \le \log(t_2) + (R+2\delta) + (k + 2R + 6\delta + 1). 
    \] 

    \medskip 
    Combining all three cases, we conclude that 
    \[ 
    d_X(x_1,x_2) \le M_{\ref{Prop_key_prop_rings}}, 
    \] 
    where \(M_{\ref{Prop_key_prop_rings}}\) is the maximum of the constants appearing above and depends on $\delta,k,R,t_1$, and $t_2$.
\end{proof}

\begin{defn}[Quasiconformal Maps, \cite{Pan1989}] \label{def_qc_maps}
    Let \((X,d_X)\) and \((Y,d_Y)\) be proper \(\delta\)-hyperbolic metric spaces for some \(\delta \ge 0\). Fix \(R > 2\delta\), and let \(\mathcal{R}(\xi,R)\) and \(\mathcal{R}(\xi',R)\) denote the ring structures on \(\partial X\) and \(\partial Y\) associated to \(\xi \in \partial X\) and \(\xi' \in \partial Y\), respectively. Let \(\eta \colon [1,\infty) \to [1,\infty)\) be a continuous function. A homeomorphism 
    \[ 
    f \colon (\partial X,\mathcal{R}(\xi,R)) \longrightarrow (\partial Y,\mathcal{R}(\xi',R)), \qquad f(\xi)=\xi', 
    \] 
     is called \(\eta\)-\emph{quasiconformal} if the following conditions are satisfied: 
    \begin{itemize}
        \item For every \(t\)-ring \(B \subset tB\) in \(\mathcal{R}(\xi,R)\), there exists an \(\eta(t)\)-ring \(B' \subset \eta(t)B'\) in \(\mathcal{R}(\xi',R)\) such that 
        \[ 
        B' \subset f(B) \subset f(tB) \subset \eta(t)B'. 
        \] 
        
        \item For every \(t\)-ring \(B' \subset tB'\) in \(\mathcal{R}(\xi',R)\), there exists an \(\eta(t)\)-ring \(B \subset \eta(t)B\) in \(\mathcal{R}(\xi,R)\) such that 
        \[ 
        B \subset f^{-1}(B') \subset f^{-1}(tB') \subset \eta(t)B. 
        \]
    \end{itemize} 
    
    A homeomorphism \(f \colon \partial X \to \partial Y\) is called \emph{quasiconformal} if there exist a constant \(R > 2\delta\), points \(\xi \in \partial X\) and \(\xi' \in \partial Y\) with \(f(\xi)=\xi'\), and a distortion function \(\eta \colon [1,\infty) \to [1,\infty)\) such that \(f \colon (\partial X,\mathcal{R}(\xi,R)) \to (\partial Y,\mathcal{R}(\xi',R)) \) is \(\eta\)-quasiconformal.
\end{defn}

Note that the composition of an \(\eta_1\)-quasiconformal map with an \(\eta_2\)-quasiconformal map is \((\eta_2 \circ \eta_1)\)-quasiconformal. The next proposition shows that quasiconformality is independent of the choice of the parameter \(R\) used to define the ring structure. In a later section, we prove that quasiconformality is also independent of the choice of the basepoint \(\xi \in \partial X\) whenever \(X\) is a visual hyperbolic space; see \autoref{prop_qc_is_independet_of_base_point}.

\begin{prop}\label{prop_ring_structures_are_independent_on_R}
    Let \(X\) be a proper \(\delta\)-hyperbolic metric space and \(\xi \in \partial X\). Given \(R, R' > 2\delta\), let \(\mathcal{R}(\xi,R)\) and \(\mathcal{R}(\xi,R')\) denote the corresponding ring structures on \(\partial X\). Then the identity map 
    \[ 
    \mathrm{id} \colon (\partial X,\mathcal{R}(\xi,R)) \longrightarrow (\partial X,\mathcal{R}(\xi,R')) 
    \] 
    is \(\eta\)-quasiconformal, where the distortion function is given by
    \(
    \eta(t) = M_{\ref{prop_ring_structures_are_independent_on_R}}t
    \) for \(t \geq 1\) and for some constant \(M_{\ref{prop_ring_structures_are_independent_on_R}} \ge 1\) depending only on \(\delta\), \(R\), and \(R'\).
\end{prop}

\begin{proof}
    Let $\mathcal{O}_\xi B(x_1, R) \subset \mathcal{O}_\xi B(x_2, R)$ be a $t$-ring in the ring structure $\mathcal{R}(\xi, R)$, centered at a point $a \in \partial X \setminus \{\xi\}$, where $t \ge 1$ and $x_1, x_2 \in [\xi,a]$ with $d_X(x_1,x_2) = \log(t)$. 
    Now, let $x_1' \in [x_1, a]$ and $x_2' \in [x_2, \xi]$ be points satisfying 
    \(
    d_X(x_1',x_1) = R' + \delta\) and \(d_X(x_2,x_2') = R + \delta.
    \)
    By \autoref{lem_Existence_of_rings}, we have
    \[
    \mathcal{O}_\xi B(x_1',R')\subset \mathcal{O}_\xi B(x_1,R) \subset \mathcal{O}_\xi B(x_2,R) \subset \mathcal{O}_\xi B(x_2',R').
    \]
    The total distance between the modified centers is given by:
    \begin{align*}
        d_X(x_1',x_2') & = d_X(x_1',x_1) + d_X(x_1,x_2) + d_X(x_2,x_2') \\
        & = (R' + \delta) + \log(t) + (R + \delta) \\
        & = \log(te^{R+R'+2\delta}).
    \end{align*}    
    
    Setting $M_{\ref{prop_ring_structures_are_independent_on_R}} = e^{R+R'+2\delta} \ge 1$, we conclude that $\mathcal{O}_\xi B(x_1',R') \subset \mathcal{O}_\xi B(x_2',R')$ is a $M_{\ref{prop_ring_structures_are_independent_on_R}}t$-ring in $\mathcal{R}(\xi,R')$.    
    In a similar way, it is follows that for any $t$-ring $\mathcal{O}_\xi B(x_1',R') \subset \mathcal{O}_\xi B(x_2',R')$ in $\mathcal{R}(\xi,R')$ there exists a $M_{\ref{prop_ring_structures_are_independent_on_R}}t$-ring $\mathcal{O}_\xi B(x_1,R) \subset \mathcal{O}_\xi B(x_2,R)$ in $\mathcal{R}(\xi,R)$ such that
    \[
    \mathcal{O}_\xi B(x_1,R)\subset \mathcal{O}_\xi B(x_1',R') \subset \mathcal{O}_\xi B(x_2',R') \subset \mathcal{O}_\xi B(x_2,R).
    \]
    Hence, the result follows. 
\end{proof}

We next recall a result of Paulin \cite{Pau96}*{Proposition 3.2}, which asserts that a quasi-isometry between proper hyperbolic metric spaces induces a quasiconformal homeomorphism between their Gromov boundaries. Alternatively, the result can also be obtained by adapting the strategy of Pansu \cite{Pan1989}*{Proposition 4.15}, who showed that a quasi-isometry between simply connected manifolds with sectional curvature satisfying \(-a^{2} < K < -b^{2} < 0\) extends homoeomorphically to the boundary at infinity as a quasiconformal map.

\begin{thm}\label{thm_qi_implies_qc}
    Given constants \(\delta,\epsilon \ge 0\), \(\lambda \ge 1\), and \(R > 2\delta\), there exists a continuous function \(\eta \colon [1,\infty) \to [1,\infty)\) such that the following holds. 
    Let \((X,d_X)\) and \((Y,d_Y)\) be proper \(\delta\)-hyperbolic metric spaces, and let \(\varphi \colon X \to Y\) be a \((\lambda,\epsilon)\)-quasi-isometry. Let $\mathcal{R}(\xi,R)$ and $\mathcal{R}(\xi', R)$ be ring structures on the Gromov boundaries $\partial X$ and $\partial Y$, respectively, where $(\xi,\xi') \in \partial X \times \partial Y$ with $\xi' = \partial \varphi (\xi)$. Then the induced boundary homeomorphism 
    \[ 
    \partial\varphi \colon (\partial X,\mathcal{R}(\xi,R)) \longrightarrow (\partial Y,\mathcal{R}(\xi',R)) 
    \] 
    is \(\eta\)-quasiconformal. Moreover, the distortion function \(\eta\) may be chosen of the form 
    \( 
    \eta(t) = C\, t^{\mu}, 
    \) for all \(t \geq 1\), 
    where the constants \(C \ge 1\) and \(\mu \ge 1\) depend only on 
    \(\delta\), \(\epsilon\), \(\lambda\), and \(R\).
\end{thm}

\section{Shadows of Horoballs} \label{sec_shadows_of_horoballs}

In this section, we study shadows of horoballs, which arise naturally from the horoball structure of the cusped space associated to a relatively hyperbolic group. We first show that the shadow of a horoball can be coarsely represented by a ring in the boundary. Building on this observation, we introduce a notion of quasiconformal maps that coarsely preserve shadows of horoballs. These shadows play a role in the study of coarsely cusp-preserving quasi-isometries of relatively hyperbolic groups and in the analysis of the quasiconformal geometry of their Bowditch boundaries.


\begin{lem}[Shadows of horoballs]\label{lem_shadow_of_horoballs}
    There exists a constant \(T_1 = T_1(\delta,R) \ge 1\) with the following property. Let \((G,\mathcal{H}_G)\) be a \(\delta\)-relatively hyperbolic group with cusped space \(X^h\), and let \(\xi \in \partial X^h\) and \(R>2\delta\). Suppose \(H \in \mathcal{H}_G\) is a horosphere-like subset with associated parabolic endpoint \(a_H \in \partial X^h \setminus \{\xi\}\).     
    Then there exists a \(T_1\)-ring \(B_H \subset T_1 B_H \) in the ring structure \(\mathcal{R}(\xi,R)\), centered at \(a_H\), such that 
    \[ 
    B_H \subset \mathcal{O}_\xi(H^h) \subset T_1 B_H . 
    \] 
    Moreover, the inner shadow \(B_H\) may be chosen of the form 
    \[ 
    B_H = \mathcal{O}_\xi B(p_{\{\xi,H\}},R), 
    \] 
    where \(p_{\{\xi,H\}} \in (\xi,a_H) \cap H^h\) is a point on a geodesic from \(\xi\) to \(a_H\) satisfying 
    \[ 
    d_X\!\left(p_{\{\xi,H\}}, H\right) \le R+\delta . 
    \] 
\end{lem}

\begin{proof}
    Let \(e_H \in H\) be the first entry point of the geodesic \((\xi,a_H)\) into the horoball \(H^h\). By \autoref{prop_visually_bounded}, there exists a constant \(K_{\ref{prop_visually_bounded}} = K_{\ref{prop_visually_bounded}}(\delta) \ge 0\) such that every geodesic ray emanating from \(\xi\) and intersecting \(H^h\) must pass through the ball \(B(e_H,K_{\ref{prop_visually_bounded}})\). It follows that 
    \[
    \mathcal{O}_\xi(H^h) \subset \mathcal{O}_\xi B(e_H,K_{\ref{prop_visually_bounded}}). 
    \]    
    
    \begin{figure}[H]
    \centering
    \begin{tikzpicture}
        \draw (0,0) circle (2.5 cm);
        \filldraw[gray!20] (0,1.7) circle (.8 cm);
        \draw (0,1.7) circle (.8 cm);

        \filldraw[gray!30] (0,0) circle (.35 cm);
        \draw (0,0) circle (.35 cm);

        \filldraw[gray!50] (0,1.3) circle (.3 cm);
        \draw (0,1.3) circle (.3 cm);

        \draw (0,-2.5)--(0,2.5);
        \draw[bend left = 10] (0,-2.5) edge (1.08,2.25);
        \draw[bend right = 10] (0,-2.5) edge (-1.08,2.25);

        \draw[-stealth] (0,.9) -- (1.2,.9);
        \draw[-stealth] (0,0) -- (1.2,0);
        \draw[-stealth] (0,1.3) -- (1.2,1.3);

        \filldraw (0,0) circle (1pt); 
        \filldraw (0,0.9) circle (1pt);
        \filldraw (0,1.3) circle (1pt);
        \node at (0,2.8) {$a_H$};
        \node at (0,-2.8) {$\xi$};
        \node at (1.8,0) {$p'_{\{\xi,H\}}$};
        \node at (1.5,.9) {$e_H$};
        \filldraw[white] (2,1.2) circle (.35cm);
        \node at (1.8,1.3) {$p_{\{\xi,H\}}$};         
    \end{tikzpicture}
    \caption{}
    \end{figure}
    
    Now choose points \(p_{\{\xi,H\}} \in [e_H,a_H )\) and \(p'_{\{\xi,H\}} \in [e_H,\xi)\) satisfying 
    \[
    d_{X^h}(p_{\{\xi,H\}},e_H) = R+\delta \quad \text{and} \quad d_{X^h}(e_H,p'_{\{\xi,H\}})=K_{\ref{prop_visually_bounded}}+\delta,
    \]
    respectively. By \autoref{lem_Existence_of_rings}, we obtain 
    \[
    \mathcal{O}_\xi B(p_{\{\xi,H\}},R) \subset \mathcal{O}_\xi (H^h) \subset \mathcal{O}_\xi B(p_{\{\xi,H\}}',R). 
    \]    
    and the distance between the centers is 
    \[
    d_{X^h}(p_{\{\xi,H\}},p'_{\{\xi,H\}}) = K_{\ref{prop_visually_bounded}} + R + 2\delta.
    \] 
    Thus, setting \(T_1 = e^{K_{\ref{prop_visually_bounded}} + R + 2\delta} \ge 1\), we conclude that \(\mathcal{O}_\xi B(p_{\{\xi,H\}}, R) \subset \mathcal{O}_\xi B(p_{\{\xi,H\}}', R)\) is a \(T_1\) ring in \(\mathcal{R}(\xi,R)\) as required.
    Finally, since \(d_{X^h}(p_{\{\xi,H\}}, e_H) = R + \delta\), the point \(p_{\{\xi, H\}}\) lies within distance \(R+\delta\) of \(H\), completing the proof.
\end{proof} 

By \autoref{lem_shadow_of_horoballs}, there exists a natural association from the collection \(\mathcal{H}_G\) of horosphere-like subsets to the ring structure \(\mathcal{R}(\xi,R)\) on \(\partial G^h\), given by 
\[
H \longmapsto \bigl( B_H \subset T_1\,B_H \bigr), 
\] 
where \(T_1 = T_1(\delta,R) \ge 1\). Moreover, by \autoref{Prop_key_prop_rings}, any two such choices of rings are uniformly comparable, with constants depending only on \(\delta\).

Now let \((G_1,\mathcal{H}_{G_1})\) and \((G_2,\mathcal{H}_{G_2})\) be \(\delta\)-relatively hyperbolic groups with Cayley graphs \(X\) and \(Y\), respectively. Fix \(R > 2\delta\) and boundary points \((\xi,\xi') \in \partial X^h \times \partial Y^h\), and let \(\mathcal{R}(\xi,R)\) and \(\mathcal{R}(\xi',R)\) denote the corresponding ring structures on \(\partial X^h\) and \(\partial Y^h\).


\begin{defn}[Coarsely Horoball-Shadow-Preserving Quasiconformal Maps]
\label{def_coarsely_shadows_of_horoballs_preserve}
    Let \(L \ge 1\). An \(\eta\)-quasiconformal homeomorphism
    \[
    f \colon (\partial X^h,\mathcal{R}(\xi,R)) \longrightarrow (\partial Y^h,\mathcal{R}(\xi',R)), \qquad f(\xi)=\xi'
    \]
    is said to \emph{\(L\)-coarsely preserve shadows of horoballs relative to \(\xi\)} if the following conditions hold:
    \begin{enumerate}
        \item Both \(f\) and \(f^{-1}\) send parabolic endpoints to parabolic endpoints.

        \item Whenever \(H \in \mathcal{H}_{G_1}\) and \(H' \in \mathcal{H}_{G_2}\) are horosphere-like subsets centered at \(a \in \partial X^h \setminus \{\xi\}\) and \(f(a) \in \partial Y^h \setminus \{\xi'\}\), respectively, satisfying
        \[
        B_H \subset \mathcal{O}_\xi(H^h) \subset T_1 B_H \quad \text{and} \quad B_{H'} \subset \mathcal{O}_{\xi'}(H'^h) \subset T_1 B_{H'},
        \]
        then
        \begin{itemize}
            \item \(\tfrac{1}{L} B_{H'} \subset f(B_H) \subset f(\mathcal{O}_\xi(H^h)) \subset f(T_1 B_H) \subset L B_{H'},\)

            \medskip

            \item \(\tfrac{1}{L} B_H \subset f^{-1}(B_{H'}) \subset f^{-1}(\mathcal{O}_{\xi'}(H'^h)) \subset f^{-1}(T_1 B_{H'}) \subset L B_H.\)
    \end{itemize}
\end{enumerate}
\end{defn}

\bigskip
\subsection{Proof of \autoref{thm_main_theorem_1}}

We recall the statement of \autoref{thm_main_theorem_1}. Let \(\varphi\to G_2\) be a coarsely cusp-preserving quasi-isometry between two relatively hyperbolic groups. Then the induced homeomorphism \(\partial\varphi^h:\partial G_1^h\to \partial G_2^h\) between their Bowditch boundaries, as described in \autoref{prop_cusp_prev_qi_implies_homeo}, is quasiconformal with respect to the associated ring structures and coarsely preserves the shadows of horoballs relative to each boundary point. 

Let \(\delta, \epsilon, K \ge 0\), \(\lambda \ge 1\), and \(R > 2\delta\). Let \((G_1,\mathcal{H}_{G_1})\) and \((G_2,\mathcal{H}_{G_2})\) be \(\delta\)-relatively hyperbolic groups with Cayley graphs \(X\) and \(Y\), respectively. Suppose that \(\varphi \colon G_1 \longrightarrow G_2\) is a \((\lambda,\epsilon)\)-quasi-isometry that is \(K\)-coarsely cusp-preserving. 
By \autoref{prop_cusp_prev_qi_implies_homeo}, there exist constants \(\lambda' \ge 1\) and \(\epsilon' \ge 0\), depending only on \(\delta,\lambda,\epsilon\), and \(K\), such that \(\varphi\) induces a \((\lambda',\epsilon')\)-quasi-isometry \( \varphi^h \colon X^h \longrightarrow Y^h, \) which is also \(K\)-coarsely cusp-preserving, together with a boundary homeomorphism \( \partial\varphi^h \colon \partial X^h \longrightarrow \partial Y^h. \)

By \autoref{thm_qi_implies_qc}, it is follows that for every choice of ring structures \(\mathcal{R}(\xi,R)\) on \(\partial X^h\) and \(\mathcal{R}(\xi',R)\) on \(\partial Y^h\), where \(\xi' = \partial\varphi^h(\xi)\), the map 
\[
\partial\varphi^h \colon
(\partial X^h,\mathcal{R}(\xi,R))
\longrightarrow
(\partial Y^h,\mathcal{R}(\xi',R))
\]
is \(\eta\)-quasiconformal for some distortion function \(\eta \colon [1,\infty) \to [1,\infty)\).

Therefore, to prove \autoref{thm_main_theorem_1}, it remains to prove that the induced boundary homeomorphism \(\partial\varphi^h\) coarsely preserves shadows of horoballs relative to each boundary point. This is established by the following proposition.

\begin{prop}\label{lem_qi_implies_coarsely_preserves_horoball_shadows}
    Let \((G_1,\mathcal{H}_{G_1})\) and \((G_2,\mathcal{H}_{G_2})\) be \(\delta\)-relatively hyperbolic groups with Cayley graphs \(X\) and \(Y\), respectively, where \(\delta \ge 0\). Let \(R>2\delta\), \(\lambda' \ge 1\), \(\epsilon' \ge 0\), and suppose that 
    \( 
    \varphi^h \colon X^h \longrightarrow Y^h 
    \)
    is a \((\lambda',\epsilon')\)-quasi-isometry. 
    Then there exists a constant \(L = L(\delta,\lambda',\epsilon',R) \ge 1\) such that, for every \(\xi \in \partial X^h\), with \(\xi' = \partial\varphi^h(\xi)\), the induced boundary homeomorphism
    \[
    \partial\varphi^h \colon (\partial X^h,\mathcal{R}(\xi,R)) \longrightarrow (\partial Y^h,\mathcal{R}(\xi',R)) 
    \] 
    \(L\)-coarsely preserves shadows of horoballs relative to \(\xi\).
\end{prop}

\begin{proof} 
    Let \(H \in \mathcal{H}_{G_1}\) be a horosphere-like subset in \(X^h\) associated to a parabolic endpoint \(a_H \in \partial X^h \setminus \{\xi\}\), and let \(H' \in \mathcal{H}_{G_2}\) be the horosphere-like subset in \(Y^h\) associated to the parabolic endpoint \(a_{H'} = \partial\varphi^h(a_H)\). By \autoref{lem_shadow_of_horoballs}, there exists a constant \(T_1 = T_1(\delta,R) \ge 1\), points \(p_{\{\xi,H\}} \in (\xi,a_H) \cap H^h\) and \(q_{\{\xi',H'\}} \in (\xi',a_{H'}) \cap H'^h\) satisfying 
    \[
    d_{X^h}(p_{\{\xi,H\}},H) \le R+\delta,
    \qquad
    d_{Y^h}(q_{\{\xi',H'\}},H') \le R+\delta,
    \]
    and \(T_1\)-rings
    \[
    B_H = \mathcal{O}_\xi B(p_{\{\xi,H\}},R) \subset T_1 B_H,
    \qquad
    B_{H'} = \mathcal{O}_{\xi'} B(q_{\{\xi',H'\}},R) \subset T_1 B_{H'},
    \]
    in \(\mathcal{R}(\xi,R)\) and \(\mathcal{R}(\xi',R)\), respectively, such that
    \[
    B_H \subset \mathcal{O}_\xi(H^h) \subset T_1 B_H,
    \qquad
    B_{H'} \subset \mathcal{O}_{\xi'}(H'^h) \subset T_1 B_{H'}.
    \]

Let \(R_0 =
K_{\ref{Prop_Stability_of_quasigeodesics}}
+ \lambda(R_1 + K_{\ref{Prop_Stability_of_quasigeodesics}})
+ 2\epsilon + 1
\). Choose a point \(x_1 \in [p_{\{\xi,H\}},a_H)\) such that
\(d_{X^h}(x_1,p_{\{\xi,H\}}) = R_0 + \delta\), and a point
\(x_2 \in [p_{\{\xi,H\}},\xi)\) such that
\(d_{X^h}(p_{\{\xi,H\}},x_2) = \log T_1\).
By \autoref{lem_Existence_of_rings}, we have
\[
\mathcal{O}_\xi B(x_1,R_0)
\subset
\mathcal{O}_\xi B(p_{\{\xi,H\}},R)=B_H
\subset
\mathcal{O}_\xi(H^h)
\subset T_1 B_H=
\mathcal{O}_\xi B(x_2,R).
\]

Using stability of quasigeodesics
(\autoref{Prop_Stability_of_quasigeodesics}), there exist points
\(y_1,y_2 \in (a_{H'},\xi')\) and a constant \(K_{\ref{Prop_Stability_of_quasigeodesics}} = K_{\ref{Prop_Stability_of_quasigeodesics}}(\delta,\lambda',\epsilon') \geq 0\) such that
\[
d_{Y^h}(\varphi^h(x_i),y_i)
\le
K_{\ref{Prop_Stability_of_quasigeodesics}},
\quad i=1,2.
\]

We claim that
\[
\mathcal{O}_{\xi'} B(y_1,R)
\subset
\partial\varphi^h\bigl(\mathcal{O}_\xi B(x_1,R_0)\bigr).
\]

Let $b' \in \mathcal{O}_{\xi'} B(y_1, R)$. There exists a geodesic $(b',\xi')$ in $Y^h$ intersecting the ball $B(y_1, R_1)$ at some point say $q \in B(y_1, R)$. Then, 
\[
d_{Y^h}(\varphi^h (x_1),q) \le R + K_{\ref{Prop_Stability_of_quasigeodesics}}.
\]
    
Suppose $(\varphi^h)^{-1} \colon Y^h \to X^h$ is a quasi-isometry inverse of $\varphi^h$. Without loss of generality, assume $(\varphi^h)^{-1}$ is also a $(\lambda,\epsilon)$-quasi-isometry and satisfies $d_X((\varphi^h)^{-1} \circ \varphi^h (x), x) \le \epsilon$, for all $x \in X^h$.    
By \autoref{Prop_Stability_of_quasigeodesics}, there exists a point $p \in (b,\xi)$, where $b = \partial (\varphi^h)^{-1}(b')$, such that $d_{X^h}(p, (\varphi^h)^{-1}(q)) \le K_{\ref{Prop_Stability_of_quasigeodesics}}$. Therefore,      
\begin{align*}
    d_{X^h}(p,x_1) & \le  d_{X^h}(p,(\varphi^h)^{-1}(q)) + d_{X^h}((\varphi^h)^{-1}(q), (\varphi^h)^{-1} \varphi^h (x_1)) + d_{X^h}((\varphi^h)^{-1} \varphi^h (x_1), x_1)\\
    & \le K_{\ref{Prop_Stability_of_quasigeodesics}} + \lambda (R_1+K_{\ref{Prop_Stability_of_quasigeodesics}}) + \epsilon + \epsilon\\
    & < R_0.
\end{align*} 
This concludes that $(b,\xi)$ intersects the ball $B(x_1,R_0)$, and so $b \in \mathcal{O}_{\xi} B(x_1,R_0)$. Hence, 
\[
\mathcal{O}_{\xi'} B(y_1,R) \subset \partial\varphi^h (\mathcal{O}_{\xi} B(x_1,R_0)).
\]
This proves the claim. Therefore, we conclude 
\[
\mathcal{O}_{\xi'} B(y_1,R) \subset \partial\varphi^h (B_H) \subset \partial\varphi^h (\mathcal{O}_{\xi}(H^h)).
\]

For the other inclusion, let \(R_1 = (2K_{\ref{Prop_Stability_of_quasigeodesics}} + \lambda R + \epsilon) + 1\) and claim that \[
\partial\varphi^h (\mathcal{O}_{\xi} B(x_2,R)) \subset  \mathcal{O}_{\xi'} B(y_2,R_1).
\]    

Let $b \in \mathcal{O}_{\xi} B(x_2,R)$, and let $(b,\xi)$ be a geodesic intersecting the ball $B(x_2,R)$ at $p$. By \autoref{Prop_Stability_of_quasigeodesics}, there exists a point $q \in (b',\xi')$, where $b' = \partial\varphi^h (b)$, such that $d_{Y^h}(\varphi^h(p),q) \le K_{\ref{Prop_Stability_of_quasigeodesics}}$. Therefore,  
\begin{align*}
    d_{Y^h}(y_2,q) & \le d_{Y^h}(y_2, \varphi^h(x_2)) + d_{Y^h}(\varphi^h (x_2), \varphi^h (p)) + d_{Y^h}(\varphi^h(p), q) \\
    & \le K_{\ref{Prop_Stability_of_quasigeodesics}} + \lambda R + \epsilon + K_{\ref{Prop_Stability_of_quasigeodesics}} \\
    & < R_1.
\end{align*}
This implies that $(b',\xi')$ intersects $B(y_2,R_1)$, and so $b' \in \mathcal{O}_{\xi'} B(y_2,R_1)$. Hence, 
\[
\partial\varphi^h(\mathcal{O}_{\xi} B(x_2,R)) \subset  \mathcal{O}_{\xi'} B(y_2,R_1).
\]
This proves the claim. Now, choose a points $y_2' \in [y_2, \xi')$ such that $d_{Y^h}(y_2, y_2') = R_1 + \delta.$ Then by \autoref{lem_Existence_of_rings}, 
\[
\mathcal{O}_{\xi'} B(y_2,R_1) \subset \mathcal{O}_{\xi'} B(y_2',R)
\]
Therefore, we conclude 
\[
 \partial\varphi^h (\mathcal{O}_{\xi}(H^h)) \subset \partial\varphi^h (T_1 B_H) \subset \mathcal{O}_{\xi'} B(y_2',R).
\]

Moreover, 
\begin{align*}
    d_{Y^h}(y_1, y_2') 
    & \leq d_{Y^h}(y_1, y_2) + d_{Y^h}(y_2, y_2') \\
    & \leq 2K_{\ref{Prop_Stability_of_quasigeodesics}} + d_{Y^h}(\varphi^h(x_1), \varphi^h(x_2)) + (R_1 + \delta) \\
    & \leq 2K_{\ref{Prop_Stability_of_quasigeodesics}} + \lambda'(R_0 + \log T_1+\delta) + \epsilon' + (R_1 + \delta).
\end{align*}

Now let \(e_H\) and \(e_{H'}\) be the first entry points of the geodesics \((a_H,\xi)\) and \((a_{H'},\xi')\) into \(H'^h\) and \(H'^h\), respectively. Since, \(R_0\) is a fix real number depending on \(\delta,\lambda,\epsilon\) and \(R\), and \(\varphi ^h\) is \(K\)-coarsely cusp-preserving \((\lambda,\epsilon)\)-quasi-isometry, then there exists a constant \(l_1 = l_1(\delta,\lambda,\epsilon,K,R) \geq 0\) such that 
\[
d_{Y^h}(y_1,e_{H'}) \le l_1.
\]

Hence, \(y_1,y_2'\) and \(q_{\{\xi',H'\}}\) lie within a uniformly bounded distance, and the bound depends only on \(\delta,\lambda',\epsilon',R\).
Choose 
\[
l_2 = \max\{d_{Y^h}(q_{\{\xi',H'\}},y_1), d_{Y^h}(q_{\{\xi',H'\}},y_2'), log(T_1)\} + R +\delta \quad \text{and} \quad L = e^{l_2}.
\] 
Then we conclude that 
\[
\tfrac{1}{L}\,B_{H'}
\subset \partial\varphi^h (B_H) \subset
\partial\varphi^h\bigl(\mathcal{O}_\xi(H^h)\bigr) \subset \partial\varphi^h (T_1 B_H)
\subset
L\,B_{H'}.
\]
Hence, \(\partial\varphi^h\) \(L\)-coarsely preserves shadows of horoballs
relative to \(\xi\). This completes the proof.

\end{proof}


\section{Quasiconformal Maps induce Quasi-Isometries} \label{sec_qc_implies_qi_for_cusped_spaces}

Through this section, we assume that \(X\) and \(Y\) are proper \(\delta\)-hyperbolic metric spaces that are \(k\)-visual with respect to a boundary point (and hence every boundary point, by \autoref{lem_visual_change_basepoint}), for some \(\delta,k\ge 0\). Fix \(R>2\delta\), and let \(\mathcal{R}(\xi,R)\) and \(\mathcal{R}(\xi',R)\) denote the ring structures on \(\partial X\) and \(\partial Y\) associated to \(\xi\in\partial X\) and \(\xi'\in\partial Y\), respectively. Suppose that
\[
f \colon (\partial X,\mathcal{R}(\xi,R))
   \longrightarrow
   (\partial Y,\mathcal{R}(\xi',R))
\]
is an \(\eta\)-quasiconformal homeomorphism satisfying \(f(\xi)=\xi'\),
where \(\eta\colon [1,\infty)\to [1,\infty)\) is a distortion function.

\begin{defn}\label{def_defining_of_F(.)}
    Assume that \(X\), \(Y\), and \(f\) are as above. Since \(X\) is \(k\)-visual with respect to \(\xi\), for every \(x\in X\) there exist a point \(z_x\) and a bi-infinite geodesic asymptotic to \(\xi\) passing through \(z_x\) such that \(d_X(x,z_x)\le k.\) We define
    \[
    F_{\{\xi,\xi'\}}(x) = \bigcup_{z_x} \left\{ y\in Y \;\middle|\; \mathcal{O}_{\xi'}B(y,R) \subset f\bigl(\mathcal{O}_{\xi}B(z_x,R)\bigr) \subset \eta(1)\,\mathcal{O}_{\xi'}B(y,R) \right\}, 
    \] 
    where the union is taken over all those points \(z_x\in B(x,k)\) each of which lies on some bi-infinite geodesic asymptotic to \(\xi\).
\end{defn}

Since \(X\) is \(k\)-visual with respect to \(\xi\) and \(f\) is
\(\eta\)-quasiconformal, it follows that \(F_{\{\xi,\xi'\}}(x)\) is a non-empty subset of \(Y\), for every \(x\in X\). With this notation and the preceding assumptions in place, the following proposition concludes that the sets
\(F_{\{\xi,\xi'\}}(x)\) have uniformly bounded diameter. 


\begin{prop} \label{prop_F(x)_is_bounded}
    For every \(r \geq 0\), there exists a constant $D_{\ref{prop_F(x)_is_bounded}} = D_{\ref{prop_F(x)_is_bounded}}(\delta,k, \eta, R,r) \ge 0$ satisfying the following property. Let $x_1,x_2 \in X$, and let \(y_1 \in F_{\{\xi,\xi'\}}(x_1)\) and \(y_2 \in F_{\{\xi,\xi'\}}(x_2)\). If \(d_{X}(x_1,x_2) \leq r\), then 
    \(
    d_{Y}(y_1,y_2) \leq D_{\ref{prop_F(x)_is_bounded}}
    \).
\end{prop} 

\begin{proof}
    Let \(y_1,\in F_{\{\xi,\xi'\}}(x_1)\) and \(y_2,\in F_{\{\xi,\xi'\}}(x_2)\). By \autoref{def_defining_of_F(.)}, for each \(i=1,2\) there exists a point \(z_{x_i} \in X\), lying on a bi-infinite geodesic asymptotic to \(\xi\), such that \(d_X(x,z_{x_i}) \le k\) and 
    \[ 
    \mathcal{O}_{\xi'} B(y_i,R) \subset f\bigl(\mathcal{O}_{\xi} B(z_{x_i},R)\bigr) \subset \eta(1)\,\mathcal{O}_{\xi'} B(y_i,R). 
    \]
    Since \(d_X(z_{x_1},z_{x_2}) \le 2k+r\), we have
    \[ 
    \mathcal{O}_{\xi} B(z_{x_1},R) \subset \mathcal{O}_{\xi} B(z_{x_2},R+2k+r) \quad \text{and} \quad \mathcal{O}_{\xi} B(z_{x_2},R) \subset \mathcal{O}_{\xi} B(z_{x_1},R+2k+r). 
    \]
    By \autoref{lem_Existence_of_rings}, setting \(t_0 = e^{R+2k+r+\delta}\), this implies 
    \[ 
    \mathcal{O}_{\xi} B(z_{x_1},R) \subset t_0\, \mathcal{O}_{\xi} B(z_{x_2},R) \quad \text{and} \quad \mathcal{O}_{\xi} B(z_{x_2},R) \subset t_0\, \mathcal{O}_{\xi} B(z_{x_1},R). 
    \]
    Since \(f\) is \(\eta\)-quasiconformal, for each \(i=1,2\) there exists an \(\eta(t_0)\)-ring 
    \[ 
    B_i' = \mathcal{O}_{\xi'} B(w_i,R) \subset \eta(t_0) B_i', \qquad w_i \in Y
    \] 
    in \(\mathcal{R}(\xi',R)\) such that 
    \[ 
    B_i' \subset f\bigl(\mathcal{O}_{\xi} B(z_{x_i},R)\bigr) \subset f\bigl(t_0 \mathcal{O}_{\xi} B(z_{x_i},R)\bigr) \subset \eta(t_0) B_i'. 
    \]
    Combining the above inclusions, we obtain: 
    \begin{itemize}
        \item \(\mathcal{O}_{\xi'} B(y_i,R) \subset \eta(t_0)\mathcal{O}_{\xi'} B(w_i,R)\) and \(\mathcal{O}_{\xi'} B(w_i,R) \subset \eta(1)\mathcal{O}_{\xi'} B(y_i,R)\) for \(i=1,2\);
        
        \item \(\mathcal{O}_{\xi'} B(w_1,R) \subset \eta(t_0)\mathcal{O}_{\xi'} B(w_2,R)\) and \(\mathcal{O}_{\xi'} B(w_2,R) \subset \eta(t_0)\mathcal{O}_{\xi'} B(w_1,R)\).
    \end{itemize}
    By \autoref{Prop_key_prop_rings}, there exist constants
    \[
    D' = D'(\delta,R,\eta(1),\eta(t_0)) \ge 0, \qquad D'' = D''(\delta,R,\eta(t_0),\eta(t_0)) \ge 0 
    \]
    such that
    \[
    d_Y(y_1,w_1) \le D', \quad d_Y(y_2,w_2) \le D', \quad d_Y(w_1,w_2) \le D''.
    \]
    Therefore,
    \[
    d_Y(y_1,y_2) \le 2D' + D''.
    \]
    By setting \(D_{\ref{prop_F(x)_is_bounded}} = 2D' + D''\), we complete the proof.
\end{proof}

Since the inverse of a $\eta$-quasiconformal homeomorphism is also \(\eta\)-quasiconformal. For each \(y\in Y\), applying \autoref{def_defining_of_F(.)} to 
\[
f^{-1}\colon (\partial Y,\mathcal{R}(\xi',R))
\longrightarrow
(\partial X,\mathcal{R}(\xi,R)),
\]
we obtain a subset
\(F_{\{\xi',\xi\}}(y)\subset X\).
As in the discussion following
\autoref{def_defining_of_F(.)}, the set \(F_{\{\xi',\xi\}}(y)\) is
nonempty for every \(y\in Y\). Moreover, applying
\autoref{prop_F(x)_is_bounded} with \(r=0\), we see that there exists a
constant \(D_{\ref{prop_F(x)_is_bounded}}=D_{\ref{prop_F(x)_is_bounded}}(\delta,k,\eta,R,0) \) such that for all \(x\in X\) and \(y\in Y\)
\[
\operatorname{diam}\bigl(F_{\{\xi,\xi'\}}(x)\bigr)\le D_{\ref{prop_F(x)_is_bounded}}
\quad\text{and}\quad
\operatorname{diam}\bigl(F_{\{\xi',\xi\}}(y)\bigr)\le D_{\ref{prop_F(x)_is_bounded}}
\]

We now define maps
\[
\Phi f\colon X\longrightarrow Y
\qquad\text{and}\qquad
\Phi f^{-1}\colon Y\longrightarrow X
\]
by assigning to each \(x\in X\) an arbitrary point of
\(F_{\{\xi,\xi'\}}(x)\), and to each \(y\in Y\) an arbitrary point of
\(F_{\{\xi',\xi\}}(y)\).
Since these sets have uniformly bounded diameter, any two choices of
\(\Phi f\) (respectively \(\Phi f^{-1}\)) differ by at most \(D_{\ref{prop_F(x)_is_bounded}}\).

With this notation and the preceding assumptions in place, the following theorem shows that
\(\Phi f\) and \(\Phi f^{-1}\) are quasi-isometries and are coarse inverses of one another.

\begin{thm}\label{thm_qc_implies_qi_for_cusped_spaces}
    There exist constants \(\lambda_1 = \lambda_1(\delta,k,\eta,R) \ge 1\) and \(\epsilon_1 = \epsilon_1(\delta,k,\eta,R) \ge 0\)
    such that the maps
    \(\Phi f \colon X \to Y \) and \(\Phi f^{-1} \colon Y \to X\) are \((\lambda_1,\epsilon_1)\)-quasi-isometries. Moreover, \(\Phi f\) and \(\Phi f^{-1}\) are coarse inverses of each other.
\end{thm}

\begin{proof}
    We first show that $\Phi f$ and $\Phi f^{-1}$ are coarse inverses of each other.

    \smallskip 
    Let $x \in X$ and set $y = \Phi f(x) \in F_{\{\xi,\xi'\}}(x)$. By \autoref{def_defining_of_F(.)}, there exists a point $z_x \in X$ lying on a bi-infinite geodesic asymptotic to $\xi$ such that
    \[
    d_X(x,z_x) \le k \quad \text{and} \quad \mathcal{O}_{\xi'} B(y,R) \subset f(\mathcal{O}_{\xi} B(z_x,R)) \subset \eta(1)\,\mathcal{O}_{\xi'} B(y,R).
    \]
    Since $f$ is $\eta$-quasiconformal, applying the same definition to $f^{-1}$, we have a point $z \in X$ such that
    \[
    \mathcal{O}_{\xi} B(z,R) \subset f^{-1}(\mathcal{O}_{\xi'} B(y,R)) \subset f^{-1}(\eta(1)\,\mathcal{O}_{\xi'} B(y,R)) \subset \eta^2(1)\,\mathcal{O}_{\xi} B(z,R).
    \]
    Consequently,
    \[
    \mathcal{O}_{\xi} B(z,R) \subset \mathcal{O}_{\xi} B(z_x,R) \subset \eta^2(1)\mathcal{O}_{\xi} B(z,R).
    \]
    By definition of $F_{\{\xi',\xi\}}(y)$, there exists $z' \in F_{\{\xi',\xi\}}(y) \subset X$ satisfying
    \[
    \mathcal{O}_{\xi} B(z',R) \subset f^{-1}(\mathcal{O}_{\xi'} B(y,R)) \subset \eta(1)\mathcal{O}_{\xi} B(z',R).
    \]
    By \autoref{Prop_key_prop_rings}, there exists a constant \(M'_{\ref{Prop_key_prop_rings}} \geq 0\) depending on \(\delta,R\) and \(\eta\) such that
    \[
    d_X(z_x,z) \le M'_{\ref{Prop_key_prop_rings}}
    \quad \text{and} \quad
    d_X(z,z') \le M'_{\ref{Prop_key_prop_rings}}.
    \]
    Moreover, by \autoref{prop_F(x)_is_bounded}, there exists a constant \(D'_{\ref{prop_F(x)_is_bounded}} = D'_{\ref{prop_F(x)_is_bounded}}(\delta,k,\eta,R,0) \geq 0\) such that
    \[
    d_X(\Phi f^{-1}(y),z') \le D'_{\ref{prop_F(x)_is_bounded}}.
    \]
    Combining these estimates, we obtain
    \begin{align*}
        d_X\bigl(x, \Phi f^{-1} \circ \Phi f(x)\bigr)
        & \le d_X(x,z_x) + d_X(z_x,z) + d_X(z,z') + d_X(z',\Phi f^{-1}(y)) \\
        & \le k + 2M'_{\ref{Prop_key_prop_rings}} + D'_{\ref{prop_F(x)_is_bounded}}.
    \end{align*}
    Setting \(\epsilon_1' = k + 2M'_{\ref{Prop_key_prop_rings}} + D'_{\ref{prop_F(x)_is_bounded}},\) we conclude that
    \[
    d_X\bigl(x, \Phi f^{-1} \circ \Phi f(x)\bigr) \le \epsilon_1'.
    \]
    A symmetric argument shows that for all $y \in Y$,
    \[
    d_Y\bigl(y, \Phi f \circ \Phi f^{-1}(y)\bigr) \le \epsilon_1'.
    \]
    Hence, $\Phi f$ and $\Phi f^{-1}$ are coarse inverses of one another.

    \smallskip 
    We now show that $\Phi f$ is a quasi-isometry.    
    Let $x_1,x_2 \in X$ with $d_X(x_1,x_2) \le 1$, and set $y_i = \Phi f(x_i)$. By \autoref{prop_F(x)_is_bounded}, there exists a constant \(D''_{\ref{prop_F(x)_is_bounded}} = D''_{\ref{prop_F(x)_is_bounded}}(\delta,k,\eta,R,1) \geq 0\) such that 
    \[
    d_Y(y_1,y_2) \le D''_{\ref{prop_F(x)_is_bounded}}.
    \]

    Now, suppose $x,z \in X$ are arbitrary. Subdivide a geodesic $[x,z]$ in $X$ from $x$ to $z$ into $x = x_0, x_1, \dots, x_{n} = z$, with $d_{X}(x_{i-1},x_i) = 1$, for $1 \le i \le (n-1)$ and $d_X(x_{n-1},x_n) \leq 1$. Then,
    \begin{align*}
        d_Y(\Phi f(x), \Phi f(z)) & \le \sum_{i=1}^n d_Y\bigl(\Phi f(x_{i-1}), \Phi f(x_{i})\bigr) \\
        & \le n\,D''_{\ref{prop_F(x)_is_bounded}} \\
        & \le D''_{\ref{prop_F(x)_is_bounded}} \bigl(d_X(x,z)+1\bigr)
    \end{align*}    
    Applying the same argument to $\Phi f^{-1}$, for all $y,w \in Y$, we have
    \[
    d_X(\Phi f^{-1}(y),\Phi f^{-1}(w)) \le D''_{\ref{prop_F(x)_is_bounded}} \bigl(d_Y(y,w)+1\bigr). 
    \]   
    
    Therefore, for any $x,z \in X$, we have
    \begin{align*}
        d_X(x,z)
        & \le d_X(x,\Phi f^{-1} \circ \Phi f(x)) + d_X(\Phi f^{-1} \circ \Phi f(x),\Phi f^{-1} \circ \Phi f(z)) \\ 
        & \qquad + d_X(\Phi f^{-1} \circ \Phi f(z),z) \\ 
        & \le 2\epsilon_1' + D''_{\ref{prop_F(x)_is_bounded}} d_Y(\Phi f(x),\Phi f(z)) + D''_{\ref{prop_F(x)_is_bounded}}. 
    \end{align*}
    Also, for all $y \in Y$, $\Phi f^{-1}(y) \in X$ satisfying $d_Y(y, \Phi f \circ \Phi f^{-1}(y)) \le \epsilon_1'$.

    Setting
    \[
    \lambda_1 = D''_{\ref{prop_F(x)_is_bounded}} \quad \text{and} \quad \epsilon_1 = 2\epsilon_1' + D''_{\ref{prop_F(x)_is_bounded}},
    \]
    we conclude that $\Phi f$ (and similarly $\Phi f^{-1}$) is a $(\lambda_1,\epsilon_1)$-quasi-isometry.
\end{proof}

Since $\Phi f \colon X \to Y$ is a quasi-isometry, by \autoref{thm_qi_implies_homeo_hyp} $\Phi f$ induces a homeomorphism 
\[
\partial \Phi f \colon \partial X \longrightarrow \partial Y
\]
between the Gromov boundaries. The following proposition shows that this induced boundary map agrees with the original quasiconformal homeomorphism \(f\). In particular, \(f\) arises as the boundary extension of a quasi-isometry between \(X\) and \(Y\).

\begin{prop}\label{prop_bdry_Pfi_same_with_f}
    The boundary map \(\partial \Phi f \colon \partial X \to \partial Y\) coincides with the given homeomorphism \(f\).
\end{prop}

\begin{proof}
    Let \(a \in \partial X \setminus \{\xi\}\), and let \((a,\xi)\) be a bi-infinite geodesic joining \(a\) to \(\xi\). To prove that \(\partial \Phi f(a) = f(a)\), it suffices to show that the image \(\Phi f((a,\xi))\) lies in a uniformly bounded neighbourhood of the geodesic \((f(a),\xi')\).
    
    Let \(x \in (a,\xi)\) be a point. By \autoref{def_defining_of_F(.)} with $z_x=x$, there exists a point \(y \in F_{\{\xi,\xi'\}}(x) \subset Y\) such that
    \[
    \mathcal{O}_{\xi'} B(y,R) \subset f(\mathcal{O}_{\xi} B(x,R)) \subset \eta(1)\,\mathcal{O}_{\xi'} B(y,R). 
    \] 
    Write \(\eta(1)\,\mathcal{O}_{\xi'} B(y,R) = \mathcal{O}_{\xi'} B(y',R)\) for some point \(y' \in Y\) lying on the same bi-infinite geodesic asymptotic to \(\xi'\) as \(y\), with \(d_Y(y,y') = \log(\eta(1)),\) and with \(y'\) closer to \(\xi'\) than \(y\). 
    Since \(a \in \mathcal{O}_{\xi} B(x,R)\), we have
    \[
    f(a) \in f(\mathcal{O}_{\xi} B(x,R)) \subset \mathcal{O}_{\xi'} B(y',R).
    \]
    Hence, the geodesic \((f(a),\xi')\) intersects the $2\delta$-neighbourhood of the ball \(B(y',R)\). This implies that the distance from \(y'\) to the geodesic \((f(a),\xi')\) is at most \(R+2\delta\). Moreover, by \autoref{prop_F(x)_is_bounded}, the distance between \(\Phi f(x)\) and \(y\) is uniformly bounded by \(D_{\ref{prop_F(x)_is_bounded}}\), for some constant \(D_{\ref{prop_F(x)_is_bounded}} = D_{\ref{prop_F(x)_is_bounded}}(\delta,k,\eta,R,0) \geq 0\). Therefore,
    \[
    d_Y\bigl(\Phi f(x),(f(a),\xi')\bigr) \le D_{\ref{prop_F(x)_is_bounded}} + \log(\eta(1)) + R + 2\delta. 
    \]

    Since \(x \in (a,\xi)\) was arbitrary, it follows that \(\Phi f((a,\xi))\) is contained in a uniformly bounded neighbourhood of the geodesic \((f(a),\xi')\). Consequently, the boundary extension \(\partial \Phi f\) sends \(a\) to \(f(a)\), and hence \(\partial \Phi f = f\).
\end{proof} 

Together, \autoref{thm_qc_implies_qi_for_cusped_spaces}, \autoref{prop_bdry_Pfi_same_with_f}, and \autoref{lem_existence_of_geodesic_in_hyp} recover Paulin's theorem for hyperbolic groups: every quasiconformal homeomorphism between the Gromov boundaries of two hyperbolic groups is arises as the boundary extension of some quasi-isometry between the groups.

Another consequence of \autoref{thm_qc_implies_qi_for_cusped_spaces}, \autoref{prop_bdry_Pfi_same_with_f}, and \autoref{thm_qi_implies_qc} is the following proposition, which shows that the notion of quasiconformality is independent of the choice of boundary points used to define the ring structures. 

\begin{prop}\label{prop_qc_is_independet_of_base_point}
    Let \(X\) and \(Y\) be proper \(\delta\)-hyperbolic metric spaces that are \(k\)-visual with respect to some (and hence every) boundary point, where \(\delta, k \ge 0\). Let \(R > 2\delta\), and let \(\mathcal{R}(\xi_0,R)\) and \(\mathcal{R}(\xi'_0,R)\) denote the ring structures on \(\partial X\) and \(\partial Y\) associated to \(\xi_0 \in \partial X\) and \(\xi'_0 \in \partial Y\), respectively. Suppose that 
    \[ 
    f \colon (\partial X,\mathcal{R}(\xi_0,R)) \longrightarrow (\partial Y,\mathcal{R}(\xi'_0,R)) 
    \] 
    is an \(\eta\)-quasiconformal homeomorphism satisfying \(f(\xi_0)=\xi'_0\), where \(\eta \colon [1,\infty) \to [1,\infty)\) is a distortion function. 
    
    Then there exists a distortion function \(\eta' \colon [1,\infty) \to [1,\infty)\), depending only on \(\delta\), \(k\), \(\eta\), and \(R\), such that for every \(\xi \in \partial X\), with \(\xi' = f(\xi)\), the homeomorphism 
    \[
    f \colon (\partial X,\mathcal{R}(\xi,R)) \longrightarrow (\partial Y,\mathcal{R}(\xi',R)) 
    \] 
    is \(\eta'\)-quasiconformal. Moreover, \(\eta'\) may be chosen of the form 
    \[ 
    \eta'(t)=Ct^\mu \qquad \text{for all } t \ge 1, 
    \] 
    where the constants \(C,\mu \ge 1\) depend only on \(\delta\), \(k\), \(\eta\), and \(R\).
\end{prop}

 In view of \autoref{lem_visual_change_basepoint}, \autoref{prop_ring_structures_are_independent_on_R}, and \autoref{prop_qc_is_independet_of_base_point}, we obtain the following corollary. It shows that the notion of quasiconformality for a boundary homeomorphism is independent of the choice of ring structures associated to corresponding boundary points used in its definition. 

\begin{cor}\label{cor_qc_is_independet_of_ring_structure}
    Let \(X\) and \(Y\) be proper \(\delta\)-hyperbolic metric spaces that are \(k\)-visual with respect to some (and hence every) boundary point, where \(\delta,k \ge 0\). Suppose that
    \[
    f \colon (\partial X,\mathcal{R}(\xi_0,R_0))
    \longrightarrow
    (\partial Y,\mathcal{R}(\xi'_0,R_0)),
    \qquad f(\xi_0)=\xi'_0,
    \]
    is an \(\eta\)-quasiconformal homeomorphism, where \(\xi_0\in\partial X\), \(\xi'_0\in\partial Y\), and \(R_0>2\delta\), and \(\eta \colon [1,\infty) \to [1,\infty)\) is a distortion function. 
    Then, for every \(R>2\delta\), there exists a distortion function \(\eta' \colon [1,\infty) \to [1,\infty)\), depending only on \(\delta\), \(k\), \(\eta\), \(R_0\), and \(R\), such that the following holds. For every \(\xi\in\partial X\), with \(\xi'=f(\xi)\), the homeomorphism \[ f \colon (\partial X,\mathcal{R}(\xi,R)) \longrightarrow (\partial Y,\mathcal{R}(\xi',R)) \] is \(\eta'\)-quasiconformal. Moreover, \(\eta'\) may be chosen of the form \[ \eta'(t)=Ct^\mu \qquad \text{for all } t\ge 1, \] where the constants \(C,\mu\ge 1\) depend only on \(\delta\), \(k\), \(\eta\), \(R_0\), and \(R\).
\end{cor}


\section{Coarsely Shadows-of-Horoballs Preserving Quasiconformal Maps Induce Coarsely Cusp-Preserving Quasi-Isometries} \label{sec_qc_implies_qi}

The following theorem shows that a quasiconformal homeomorphism between Bowditch boundaries, which coarsely preserves shadows of horoballs, with respect to every corresponding pair of boundary points, induces a coarsely cusp-preserving quasi-isometry between the associated cusped spaces. 

\begin{thm}\label{thm_cusp_preserving_qi_btwn_cusped_spaces} 
    Let \((G_1,\mathcal H_{G_1})\) and \((G_2,\mathcal H_{G_2})\) be \(\delta\)-relatively hyperbolic groups with Cayley graphs \(X\) and \(Y\), respectively, and let \(X^h\) and \(Y^h\) be the corresponding cusped spaces. 
    Fix constants \(R>2\delta\) and \(L\ge 1\) and \(L\ge 1\), and let \(f\colon \partial X^h \to \partial Y^h\) be a homeomorphism. Suppose that for every \(\xi\in\partial X^h\), with \(\xi'=f(\xi)\), the homeomorphism 
    \[ 
    f\colon (\partial X^h,\mathcal R(\xi,R)) \longrightarrow (\partial Y^h,\mathcal R(\xi',R)) 
    \] 
    is \(\eta\)-quasiconformal and \(L\)-coarsely preserves shadows of horoballs relative to \(\xi\), where \( \eta\colon [1,\infty)\to[1,\infty) \) is a distortion function. Then there exists a constant \( K_{\ref{thm_cusp_preserving_qi_btwn_cusped_spaces}} = K_{\ref{thm_cusp_preserving_qi_btwn_cusped_spaces}} (\delta,\eta,R,L)\ge 0 \) such that \(f\) induces a quasi-isometry \[ \Phi f\colon X^h\longrightarrow Y^h \] which \(K_{\ref{thm_cusp_preserving_qi_btwn_cusped_spaces}}\)-coarsely preserves cusps. 
\end{thm} 

Before proving \autoref{thm_cusp_preserving_qi_btwn_cusped_spaces}, we establish the following lemma, which will be needed. 
The lemma concludes that, although the quasi-isometries that are constructed using different choices of boundary points, they remain at uniformly bounded distance from one another.

\begin{lem} \label{lem_two_qi_with_same_bdry_are_bdd}
    Let \(\delta \ge 0\), \(\lambda \ge 1\), and \(\epsilon \ge 0\). Let \((G_1,\mathcal{H}_{G_1})\) and \((G_2,\mathcal{H}_{G_2})\) be two \(\delta\)-relatively hyperbolic groups with Cayley graphs \(X\) and \(Y\), respectively. Suppose that
    \[
    \phi_1 \colon X^h \longrightarrow Y^h
    \quad \text{and} \quad
    \phi_2 \colon X^h \longrightarrow Y^h
    \]
    are \((\lambda,\epsilon)\)-quasi-isometries, respectively, such that their boundary extension maps on the Bowditch boundary coincide, that is,
    \(
    \partial \phi_1 = \partial \phi_2.
    \)
    Then there exists a constant \(M_{\ref{lem_two_qi_with_same_bdry_are_bdd}} = M_{\ref{lem_two_qi_with_same_bdry_are_bdd}}(\delta,\lambda,\epsilon) \geq 0\) such that
    \[
    \sup_{x \in X^h} d_{Y^h}\bigl(\phi_1(x),\phi_2(x)\bigr)
    \le
    M_{\ref{lem_two_qi_with_same_bdry_are_bdd}}.
    \]
\end{lem}

To prove \autoref{lem_two_qi_with_same_bdry_are_bdd}, we use the notion of quasi-centres together with some of their basic properties, which we briefly recall.
Let \(X\) be a proper \(\delta\)-hyperbolic metric space, with \(\delta \ge 0\). Given a triple \((a,b,c)\) of pairwise distinct points in \(X \cup \partial X\), there exists a point \(p_{abc} \in X\), called a \emph{quasi-centre}, such that the distance from \(p_{abc}\) to each side of any geodesic triangle \(\triangle(a,b,c)\) is at most \(5\delta\). Moreover, any two quasi-centres associated to the same triple \((a,b,c)\) are at uniformly bounded distance from each other, with a bound depending only on \(\delta\).

Quasi-isometries between hyperbolic metric spaces coarsely preserve quasi-centres; see \cite{Mackay-Sisto}*{Lemma~3.14}. Furthermore, if \(X\) is the Cayley graph of a hyperbolic group, or the cusped space of a relatively hyperbolic group, then for every point \(x \in X\) there exists a triple \((a,b,c)\) of pairwise distinct points in \(\partial X \) such that \(x\) lies within a uniformly bounded distance of a quasi-centre associated to \((a,b,c)\). This bound depends only on the hyperbolicity constant. For further details on quasi-centres and their properties, see \cites{Pau96,Mackay-Sisto,PS2024}. 

\begin{proof}[Proof of \autoref{lem_two_qi_with_same_bdry_are_bdd}]
    Let \(x \in X^h\). By the discussion above, there exist three pairwise distinct points \(a,b,c \in \partial X^h\) such that \(x\) lies within a uniformly bounded distance of a quasi-centre associated to the triple \((a,b,c)\), where the bound depends only on \(\delta\).

    Let \(f = \partial \phi_1 = \partial \phi_2\), and let \(q \in Y^h\) be a quasi-centre corresponding to the triple \((f(a),f(b),f(c))\). Since quasi-isometries coarsely preserve quasi-centres, each point \(\phi_i(x)\), for \(i = 1,2\), lies within a uniformly bounded distance of \(q\), where the bound depends only on \(\delta\), \(\lambda\), and \(\epsilon\). Consequently, the distance between \(\phi_1(x)\) and \(\phi_2(x)\) is uniformly bounded by a constant depending only on \(\delta\), \(\lambda\), and \(\epsilon\). Since \(x \in X^h\) was arbitrary, the conclusion of the lemma follows.
\end{proof}


\medskip

\begin{proof}[Proof of \autoref{thm_cusp_preserving_qi_btwn_cusped_spaces}]
    For every pair \((\xi,\xi') \in \partial X^h \times \partial Y^h\) with \(\xi' = f(\xi)\), since the map
    \[
    f \colon (\partial X^h,\mathcal{R}(\xi,R)) \longrightarrow (\partial Y^h,\mathcal{R}(\xi',R))
    \]
    is \(\eta\)-quasiconformal, by \autoref{lem_existence_of_geodesic_in_rel_hyp}, \autoref{thm_qc_implies_qi_for_cusped_spaces} and \autoref{prop_bdry_Pfi_same_with_f}, there exist constants \(\lambda_1 \ge 1\) and \(\epsilon_1 \ge 0\) depend only on \(\delta,\eta\), and \(R\) such that the map \(f\) induces a \((\lambda_1,\epsilon_1)\)-quasi-isometry
    \[
    \Phi f_{\{\xi,\xi'\}} \colon X^h \longrightarrow Y^h \quad \text{satisfying} \quad \partial\Phi f_{\{\xi,\xi'\}} = f.
    \]
    Although these quasi-isometries are constructed using different basepoints, their extended boundary homeomorphisms coincide with the given quasiconformal map $f$, and hence, by \autoref{lem_two_qi_with_same_bdry_are_bdd} they differ from one another by at most a uniformly bounded distance $M_{\ref{lem_two_qi_with_same_bdry_are_bdd}}$ depending on $\delta,\lambda$, and $\epsilon$. We fix a pair \((\xi_0,\xi_0') \in \partial X^h \times \partial Y^h\) with \(\xi_0' = f(\xi_0)\), and define
    \[
    \Phi f = \Phi f_{\{\xi_0,\xi_0'\}}.
    \]

    Then \(\Phi f\) is a \((\lambda_1,\epsilon_1)\)-quasi-isometry, and by \autoref{prop_bdry_Pfi_same_with_f}, the boundary extension map satisfy $\partial \Phi f = f$. 
    By \autoref{lem_two_qi_with_same_bdry_are_bdd}, there exists a constant
    \( M_{\ref{lem_two_qi_with_same_bdry_are_bdd}} = M_{\ref{lem_two_qi_with_same_bdry_are_bdd}}(\delta,\eta,R) \ge 0
    \)
    such that for every \((\xi,\xi')\) with \(\xi' = f(\xi)\),
    \[
    d_{Y^h}\bigl(\Phi f(x), \Phi_{\{\xi,\xi'\}}(x)\bigr)
    \le M_{\ref{lem_two_qi_with_same_bdry_are_bdd}}
    \quad \text{for all } x \in X^h.
    \]

    We now show that \(\Phi f\) coarsely preserves cusps.
    Let \(H \in \mathcal{H}_{G_1}\) be a horosphere-like subset of \(X^h\) with parabolic endpoint \(a_H \in \partial X^h\), and let \(x \in H\).
    There exists a bi-infinite geodesic \((a_H,\xi)\) passing through \(x\), where \(\xi \in \partial X^h \setminus \{a_H\}\).
    By \autoref{lem_shadow_of_horoballs}, there exists a constant \(T_1 = T_1(\delta,R) \ge 1\) and a \(T_1\)-ring
    \(
    B_{\{\xi,H\}} = \mathcal{O}_{\xi} B(p_{\{\xi,H\}},R) \subset T_1 B_{\{\xi,H\}}
    \)
    in \(\mathcal{R}(\xi,R)\) with \(p_{\{\xi,H\}} \in (a_H,\xi) \cap H^h\) and with \(d_{X^h}(p_{\{\xi,H\}},H) \le R + \delta\) such that
    \[
    B_{\{\xi,H\}} \subset \mathcal{O}_{\xi}(H^h) \subset T_1 B_{\{\xi,H\}}.
    \]
    
    Let \(H' \in \mathcal{H}_{G_2}\) be the horosphere-like subset of \(Y^h\) centered at
    \(a_{H'} = f(a_H)\).
    Again by \autoref{lem_shadow_of_horoballs}, there exists a \(T_1\)-ring
    \(
    B_{\{\xi',H'\}} = \mathcal{O}_{\xi'} B(q_{\{\xi',H'\}},R) \subset T_1 B_{\{\xi',H'\}}
    \)
    in \(\mathcal{R}(\xi',R)\) satisfying
    \[
    B_{\{\xi',H'\}} \subset \mathcal{O}_{\xi'}(H'^h) \subset T_1 B_{\{\xi',H'\}}.
    \]
    Since \(f\) \(L\)-coarsely preserves shadows of horoballs relative to \(\xi\), we obtain
    \[
    \tfrac{1}{L} B_{\{\xi',H'\}} \subset f(B_{\{\xi,H\}})
    \subset f(\mathcal{O}_{\xi}(H^h)) \subset f(T_1 B_{\{\xi,H\}})
    \subset L B_{\{\xi',H'\}}.
    \]
    Let \(\tfrac{1}{L} B_{\{\xi',H'\}} = \mathcal{O}_{\xi'} B(q,R)\), where
    \(q \in [q_{\{\xi',H'\}},a_{H'})\) and \(d_{Y^h}(q,q_{\{\xi',H'\}}) = \log L\).
    Then
    \[
    \mathcal{O}_{\xi'} B(q,R)
    \subset f(\mathcal{O}_{\xi} B(p_{\{\xi,H\}},R))
    \subset L^2 \mathcal{O}_{\xi'} B(q,R).
    \]
    By the definition of \(F_{\{\xi,\xi'\}}\) (\autoref{def_defining_of_F(.)}), there exists
    \(y \in F_{\{\xi,\xi'\}}(p_{\{\xi,H\}})\) such that
    \[
    \mathcal{O}_{\xi'} B(y,R)
    \subset f(\mathcal{O}_{\xi} B(p_{\{\xi,H\}},R))
    \subset \eta(1)\mathcal{O}_{\xi'} B(y,R).
    \]
    Applying \autoref{Prop_key_prop_rings}, there exists a constant
    \(M_{\ref{Prop_key_prop_rings}} = M_{\ref{Prop_key_prop_rings}}(\delta,R,L^2,\eta(1))\)
    such that
    \[
    d_{Y^h}(y,q)
    \le M_{\ref{Prop_key_prop_rings}},
    \]
    which implies
    \[
    d_{Y^h}(y,q_{\{\xi',H'\}})
    \le M_{\ref{Prop_key_prop_rings}} + \log L.
    \]

    By \autoref{lem_existence_of_geodesic_in_rel_hyp}, \(Y^h\) is \(3\delta\)-visual with respect to \(\xi'\), and hence by \autoref{prop_F(x)_is_bounded}, there exists
    \(D_{\ref{prop_F(x)_is_bounded}} \ge 0\) depending on \(\delta,\eta\) and \(R\) such that
    \[
    d_{Y^h}\bigl(y,\Phi_{\{\xi,\xi'\}}(p_{\{\xi,H\}})\bigr)
    \le D_{\ref{prop_F(x)_is_bounded}}.
    \]

    Combining these estimates,    
    \begin{align*}
        d_{Y^h}(\Phi f_{\{\xi,\xi'\}} (x), H') 
        &\le d_{Y^h}(\Phi f_{\{\xi,\xi'\}} (x), \Phi f_{\{\xi,\xi'\}}(p_{\{\xi,H\}})) \\
        & \qquad + d_{Y^h}(\Phi f_{\{\xi,\xi'\}}(p_{\{\xi,H\}}), q_{\{\xi',H'\}}) \\
        & \qquad + d_{Y^h}(q_{\{\xi',H'\}}, H') \\
        & \le (\lambda_1(R + \delta) + \epsilon_1) + (D_{\ref{prop_F(x)_is_bounded}} + M_{\ref{Prop_key_prop_rings}} + \log L) + (R+\delta).
    \end{align*}
    
    Finally, since \(\Phi f\) and \(\Phi_{\{\xi,\xi'\}}\) differ by at most
    \(M_{\ref{lem_two_qi_with_same_bdry_are_bdd}}\),
    we obtain
    \[
    d_{Y^h}(\Phi f(x),H')
    \le K'_{\ref{thm_cusp_preserving_qi_btwn_cusped_spaces}} + M_{\ref{lem_two_qi_with_same_bdry_are_bdd}},
    \]
    where \(K_{\ref{thm_cusp_preserving_qi_btwn_cusped_spaces}}' = (\lambda_1(R + \delta) + \epsilon_1) + (D_{\ref{prop_F(x)_is_bounded}} + M_{\ref{Prop_key_prop_rings}} + \log(L)) + (R+\delta)\) depends only on
    \(\delta,\eta,R\), and \(L\). Thus, we conclude 
    \[
    \Phi f(H) \subset N_{K_{\ref{thm_cusp_preserving_qi_btwn_cusped_spaces}}' + M_{\ref{lem_two_qi_with_same_bdry_are_bdd}}} (H').
    \]

    A symmetric argument applied to the inverse quasi-isometry $\Phi f^{-1}: Y^h \to X^h$, we conclude that for all $H \in \mathcal{H}_{G_2}$, there exists $H' \in \mathcal{H}_{G_1}$ such that 
    \[
    \Phi f^{-1}(H) \subset N_{K_{\ref{thm_cusp_preserving_qi_btwn_cusped_spaces}}' + M_{\ref{lem_two_qi_with_same_bdry_are_bdd}}} (H').
    \]
    Since $\Phi f$ and $\Phi f^{-1}$ are $\epsilon_1$-coarse inverse of each other, the preimage of $H$ under $\Phi f$ satisfies
    \[
    (\Phi f)^{-1}(H) \subset N_{K_{\ref{thm_cusp_preserving_qi_btwn_cusped_spaces}}' + M_{\ref{lem_two_qi_with_same_bdry_are_bdd}} + \epsilon_1} (H').
    \]

    Therefore, setting $K_{\ref{thm_cusp_preserving_qi_btwn_cusped_spaces}} = K_{\ref{thm_cusp_preserving_qi_btwn_cusped_spaces}}' + M_{\ref{lem_two_qi_with_same_bdry_are_bdd}} + \epsilon_3$, we conclude $\Phi f$ (and similarly $\Phi f^{-1}$) are $K_{\ref{thm_cusp_preserving_qi_btwn_cusped_spaces}}$-coarsely cusp-preserving quasi-isometry.
\end{proof}

The quasi-isometries \(\Phi f\) and \(\Phi f^{-1}\) naturally induce maps
\(\varphi_f \colon G_1 \to G_2\) and \(\varphi_{f^{-1}} \colon G_2 \to G_1\), respectively, as we now describe.
Since \(\Phi f\) (respectively \(\Phi f^{-1}\)) is
\(K_{\ref{thm_cusp_preserving_qi_btwn_cusped_spaces}}\)-coarsely cusp-preserving,
the image \(\Phi f(G_1)\) (respectively \(\Phi f^{-1}(G_2)\)) is contained in the
\(K_{\ref{thm_cusp_preserving_qi_btwn_cusped_spaces}}\)-neighbourhood of \(G_2\)
(respectively \(G_1\)).
We define
\[
\varphi_f \coloneqq \operatorname{Pr}_{G_2} \circ \Phi f \restriction_{G_1}
\quad \text{and} \quad
\varphi_{f^{-1}} \coloneqq \operatorname{Pr}_{G_1} \circ \Phi f^{-1} \restriction_{G_2},
\]
where, for \(i = 1,2\),
\(
\operatorname{Pr}_{G_i} \colon
N_{K_{\ref{thm_cusp_preserving_qi_btwn_cusped_spaces}}}(G_i) \to G_i
\)
denotes a nearest-point projection map defined on the
\(K_{\ref{thm_cusp_preserving_qi_btwn_cusped_spaces}}\)-neighbourhood
\(N_{K_{\ref{thm_cusp_preserving_qi_btwn_cusped_spaces}}}(G_i)\).
We now show that \(\varphi_f\) is the desired coarsely cusp-preserving
quasi-isometry satisfying the conclusion of
\autoref{thm_main_theorem_2}.

\begin{thm}\label{thm_cusp_preserving_qi_implies_ambient_qi_for_QC}
    There exist constants \(\epsilon, K \ge 0\) and \(\lambda \ge 1\), depending only on \(\delta, \lambda_1, \epsilon_1\), and \(K_{\ref{thm_cusp_preserving_qi_btwn_cusped_spaces}}\) (equivalently, only on \(\delta, \eta\), and \(R\)), such that the map 
    \[
    \varphi_f \colon (G_1, d_{G_1}) \longrightarrow (G_2, d_{G_2})
    \]
    is a \(K\)-coarsely cusp-preserving \((\lambda, \epsilon)\)-quasi-isometry. Moreover, \(\varphi_{f^{-1}}\) is a quasi-isometry inverse of \(\varphi_f\). 
\end{thm}

\begin{proof}   
    From the above construction of $\varphi_f$ and $\varphi_{f^{-1}}$, it follows that for every $x\in G_1$ and $y \in G_2$ 
    \[
    d_{Y^h}(\Phi f(x), \varphi_f(x)) \le K_{\ref{thm_cusp_preserving_qi_btwn_cusped_spaces}} \quad \text{and} \quad d_{X^h}(\Phi f^{-1}(y), \varphi_{f^{-1}}(y)) \le K_{\ref{thm_cusp_preserving_qi_btwn_cusped_spaces}}.
    \] 


    \noindent\textbf{Claim 1.} First, we show that $\varphi_f$ and $\varphi_{f^{-1}}$ are coarse inverses of each other. 
    
    Let $x \in G_1$. Then,
    \begin{align*}
        & d_{Y^h}(\Phi f (x), \varphi_f (x)) \le K_{\ref{thm_cusp_preserving_qi_btwn_cusped_spaces}} \\
        \implies & d_{X^h}(\Phi f^{-1} \circ \Phi f (x), \Phi f^{-1} \circ \varphi_f (x)) \le \lambda_1 K_{\ref{thm_cusp_preserving_qi_btwn_cusped_spaces}} + \epsilon_1 \\
        \implies & d_{X^h}(x, \Phi f^{-1} \circ \varphi_f (x)) \le (\lambda_1 K_{\ref{thm_cusp_preserving_qi_btwn_cusped_spaces}} + \epsilon_1) + \epsilon_1 \\
        \implies & d_{X^h}(x, \varphi_{f^{-1}} \circ \varphi_f (x)) \le (\lambda_1 K_{\ref{thm_cusp_preserving_qi_btwn_cusped_spaces}} + 2\epsilon_1) + K_{\ref{thm_cusp_preserving_qi_btwn_cusped_spaces}}
    \end{align*}
    Since $X$, and hence $G_1$, is properly embedded in $X^h$, by \autoref{lem_Uniformly_Properly_embedded_Lemma}, there exists a constant $K_{\ref{lem_Uniformly_Properly_embedded_Lemma}}' = K_{\ref{lem_Uniformly_Properly_embedded_Lemma}}'(\delta, \lambda_1 K_{\ref{thm_cusp_preserving_qi_btwn_cusped_spaces}} + 2\epsilon_1 + K_{\ref{thm_cusp_preserving_qi_btwn_cusped_spaces}}) \ge 0$ such that for every $x \in G_1$,
    \[
    d_{G_1}(x, \varphi_{f^{-1}} \circ \varphi_f (x)) \le K_{\ref{lem_Uniformly_Properly_embedded_Lemma}}'.
    \]
    In a similar way, we can show that for every $y \in G_2$ 
    \[
    d_{G_2}(y, \varphi_f \circ \varphi_{f^{-1}} (y)) \le K_{\ref{lem_Uniformly_Properly_embedded_Lemma}}'.
    \]
    Hence, $\varphi_f$ and $\varphi_{f^{-1}}$ are coarse inverses of each other.


    \medskip
    \noindent\textbf{Claim 2.} We now show that $\varphi_f$ is a quasi-isometry. 
    
    Let $x,z \in G_1$ satisfy $d_{G_1}(x,z) = 1$. Then there exist $H_1, H_2 \in \mathcal{H}_{G_1}$ such that $x \in H_1$ and $z \in H_2$. Since $\Phi f$ is a $(\lambda_1,\epsilon_1)$-quasi-isometry, we have
    \[
    d_{Y^h}(\Phi f (x), \Phi f (z)) \le \lambda_1 + \epsilon_1 \quad \implies \quad d_{Y^h}(\varphi_f (x), \varphi_f (z)) \le \lambda_1 + \epsilon_1 + 2K_{\ref{thm_cusp_preserving_qi_btwn_cusped_spaces}}.
    \]
    Since $Y$, and hence $G_2$, is properly embedded in $Y^h$, it follows from \autoref{lem_Uniformly_Properly_embedded_Lemma} that there exists a constant $K_{\ref{lem_Uniformly_Properly_embedded_Lemma}}'' = K_{\ref{lem_Uniformly_Properly_embedded_Lemma}}''(\delta, \lambda_1 + \epsilon_1 + 2K_{\ref{thm_cusp_preserving_qi_btwn_cusped_spaces}}) \ge 0$ such that 
    \[
    d_{G_2}(\varphi_f (x), \varphi_f (z)) \le K_{\ref{lem_Uniformly_Properly_embedded_Lemma}}''.
    \]

    Now, let $x,z \in G_1$ be arbitrary. Partition the geodesic segment $[x,z]_{G_1}$ in $G_1$ into $x=x_0,x_1,\dots,x_{n}=z$, with $d_{G_1}(x_{i-1},x_i) = 1$  for $1 \le i \le n$. Then,
    \begin{align*}
        d_{G_2}(\varphi_f (x), \varphi_f (z)) & \le \sum_{i=1}^{n} d_{G_2}(\varphi_f (x_{i-1}), \varphi_f (x_i)) \\
        & \le n K_{\ref{lem_Uniformly_Properly_embedded_Lemma}}'' \\
        & = K_{\ref{lem_Uniformly_Properly_embedded_Lemma}}'' d_{G_1}(x,z)
    \end{align*}
    
    Similarly, using the quasi-isometry inverse  $\Phi f^{-1}$, for all $y,w \in G_2$, we have 
    \[
    d_{G_1}(\varphi_{f^{-1}} y, \varphi_{f^{-1}} w) \le K_{\ref{lem_Uniformly_Properly_embedded_Lemma}}'' d_{G_2}(y,w).
    \]
    Therefore, for all $x,z \in G_1$
    \begin{align*}
        & d_{G_1}(\varphi_{f^{-1}} \circ \varphi_f(x), \varphi_{f^{-1}} \circ \varphi_f(z)) \le K_{\ref{lem_Uniformly_Properly_embedded_Lemma}}'' d_{G_2}(\varphi_f(x),\varphi_f(z))\\
        \implies &  d_{G_1}(x,z) \le K_{\ref{lem_Uniformly_Properly_embedded_Lemma}}'' d_{G_2}(\varphi_f(x),\varphi_f(z)) + 2 K_{\ref{lem_Uniformly_Properly_embedded_Lemma}}'.
    \end{align*}

    Also, for all $y \in G_2$, we have 
    \(
    \varphi_{f^{-1}}(y) \in G_1\) and \(d_{G_2}(y,\varphi_f \circ \varphi_{f^{-1}}(y)) \le K_{\ref{lem_Uniformly_Properly_embedded_Lemma}}'.
    \)
    Hence, setting $\lambda = K_{\ref{lem_Uniformly_Properly_embedded_Lemma}}''$  and $\epsilon = 2K_{\ref{lem_Uniformly_Properly_embedded_Lemma}}'$, we conclude that $\varphi_f$ (similarly $\varphi_{f^{-1}}$) is a $(\lambda,\epsilon)$-quasi-isometry. 
    
    \medskip
    \noindent\textbf{Claim 3.} Finally, we show that $\varphi_f$ is a coarsely cusp-preserving map.

    For every $H \in \mathcal{H}_{G_1}$, there exists $H' \in \mathcal{H}_{G_2}$ such that $\Phi (H) \subset N_{K_{\ref{thm_cusp_preserving_qi_btwn_cusped_spaces}}}(H')$, which implies 
    \[
    \varphi_f (H) \subset N_{2K_{\ref{thm_cusp_preserving_qi_btwn_cusped_spaces}}}(H').
    \]
    Conversely, for all $H \in \mathcal{H}_{G_2}$, there exists $H' \in \mathcal{H}_{G_1}$ such that
    \begin{align*}
        \varphi_f (x) \in H & \implies \Phi f (x) \in N_{K_{\ref{thm_cusp_preserving_qi_btwn_cusped_spaces}}}(H) \\ 
        & \implies \Phi f^{-1} \circ \Phi f(x) \in N_{ \lambda_1 K_{\ref{thm_cusp_preserving_qi_btwn_cusped_spaces}} + \epsilon_1 + K_{\ref{thm_cusp_preserving_qi_btwn_cusped_spaces}}} (H') \\
        & \implies x \in N_{ \lambda_1 K_{\ref{thm_cusp_preserving_qi_btwn_cusped_spaces}} + 2\epsilon_1 + K_{\ref{thm_cusp_preserving_qi_btwn_cusped_spaces}}} (H') 
    \end{align*}
    Therefore, 
    \[
    (\varphi_f)^{-1}(H) \subset N_{ \lambda_1 K_{\ref{thm_cusp_preserving_qi_btwn_cusped_spaces}} + 2\epsilon_1 + K_{\ref{thm_cusp_preserving_qi_btwn_cusped_spaces}}} (H').
    \]
    Hence, setting $K = \lambda_1 K_{\ref{thm_cusp_preserving_qi_btwn_cusped_spaces}} + 2\epsilon_1 + K_{\ref{thm_cusp_preserving_qi_btwn_cusped_spaces}}$, we conclude that $\varphi_f$ (and similarly $\varphi_{f^{-1}}$) preserves cusps $K$-coarsely. 
    This proves the claim.

    The required conclusion now follows by combining Claims~1, 2, and~3.
\end{proof} 

\begin{prop}\label{prop_boundaries_coincide}
    Let \(f:\partial X^h\to\partial Y^h\) and $\varphi_f: X \to Y$ are as in \autoref{thm_cusp_preserving_qi_implies_ambient_qi_for_QC}. Let $\varphi_f^h: X^h \to Y^h$ denote its extension as defined in \autoref{prop_cusp_prev_qi_implies_homeo}. Then \(\partial\varphi_f^h = f\). 
\end{prop}

\begin{proof}
    By \autoref{prop_cusp_prev_qi_implies_homeo}, the map \(\varphi_f^h \colon X^h \to Y^h\) is a coarsely cusp-preserving quasi-isometry. Let \(\gamma \colon [0,\infty) \to X^h\) be a vertical geodesic ray in $X^h$ with \(\gamma(0) \in H\). Then the image \(\varphi_f^h(\gamma)\) is eventually a vertical geodesic ray, while $\Phi f(\gamma)$ is a quasigeodesic ray. Moreover, these two rays remain within a uniformly bounded Hausdorff distance of one another. It follows that for every \(x \in X^h\), the distance in \(Y^h\) between \(\varphi_f^h(x)\) and \(\Phi f(x)\) is uniformly bounded. Consequently, the induced boundary map \(\partial \varphi_f^h \colon \partial X^h \to \partial Y^h\), extension of \(\varphi_f^h\), coincides with \(\partial \Phi f\), and hence with the given homeomorphism \(f\).
\end{proof}


\bibliography{Bibliography.bib}


\end{document}